\documentclass[10pt]{article}
\usepackage[latin1]{inputenc}
\usepackage{epsfig}
\usepackage{color}
\usepackage[british,english]{babel}
\usepackage{amsthm}
\usepackage{amsmath}
\usepackage{amsfonts}
\usepackage{amssymb}
\usepackage{graphicx}

\usepackage{fancyhdr}

\newtheorem{corollary}{Corollary}[section]
\newtheorem{definition}[corollary]{Definition}

\newtheorem{lemma}[corollary]{Lemma}
\newtheorem{proposition}[corollary]{Proposition}
\newtheorem{remark}[corollary]{Remark}
\newtheorem{theorem}[corollary]{Theorem}

\date{}
\begin{document}
\title{Homogenization of parabolic problems with dynamical boundary conditions of reactive-diffusive type in perforated media}\maketitle

\vskip-30pt
 \centerline{Mar\'ia ANGUIANO\footnote{Departamento de An\'alisis Matem\'atico. Facultad de Matem\'aticas.
Universidad de Sevilla, 41012 Sevilla (Spain)
anguiano@us.es}}

 \renewcommand{\abstractname} {\bf Abstract}
\begin{abstract} 
This paper deals with the homogenization of the reaction-diffusion equations in a domain containing periodically distributed holes of size $\varepsilon$, with a dynamical boundary condition of {\it reactive-diffusive} type, i.e., we consider the following nonlinear boundary condition on the surface of the holes 
$$
\nabla u_\varepsilon \cdot \nu+\varepsilon\,\displaystyle\frac{\partial u_\varepsilon}{\partial
t}=\varepsilon\,\delta \Delta_{\Gamma}u_\varepsilon-\varepsilon\,g(u_\varepsilon),
$$
where $\Delta_{\Gamma}$ denotes the Laplace-Beltrami operator on the surface of the holes, $\nu$ is the outward normal to the boundary, $\delta>0$ plays the role of a surface diffusion coefficient and $g$ is the nonlinear term. We generalize our previous results (see \cite{Anguiano}) established in the case of a dynamical boundary condition of {\it pure-reactive} type, i.e., with $\delta=0$. We prove the convergence of the homogenization process to a nonlinear reaction-diffusion equation whose diffusion matrix takes into account the reactive-diffusive condition on the surface of the holes. 
\end{abstract}

 {\small \bf AMS classification numbers: 35B27, 35K57}  
 
 
 \section {Introduction and setting of the problem}\label{S1} 
In the context of reaction-diffusion equations, dynamical boundary conditions have been rigorously derived in Gal and Shomberg \cite{Gal2} based on first and second thermodynamical principles and their physical interpretation was also given in Goldstein \cite{Goldstein}. It is worth emphasizing that the derivation in \cite{Gal2} obtains the dynamical boundary condition of {\it reactive-diffusive} type both as a sufficient and necessary condition for thermodynamic processes which incorporate thermodynamic sources located along the boundary, and in which the second law plays a crucial role, while in \cite{Goldstein} it has been introduced only as a sufficient condition. 

In particular, a dynamical boundary condition of {\it reactive-diffusive} type accounts for (see \cite[Section 3]{Goldstein}) a heat source on the boundary that can depend on the heat flow along the boundary, the heat flux across the boundary and the temperature at the boundary.
Consider the reaction-diffusion equation, with dynamical boundary condition of {\it reactive-diffusive} type, provides, in addition to classical bulk diffusion, a diffusion mechanism present along the boundary. A typical example in the theory of heat conduction (see \cite{Gal2} for more details) arises when a given body is in perfect thermal contact with a sufficiently thin metal sheet, possibly of different material and completely insulating the internal body from external contact, say, a well-stirred hot or cold fluid. 
 
In a recent article (see \cite{Anguiano}) we addressed the problem of the homogenization of the reaction-diffusion equations with a dynamical boundary condition of {\it pure-reactive} type in a domain perforated with holes. The present article is devoted to the generalization of that previous study to the case of a dynamical boundary condition of {\it reactive-diffusive} type, i.e., we add to the dynamical boundary condition a Laplace-Beltrami correction term. Let us introduce the model we will be involved with in this article. 
 
{\bf The geometrical setting.} Let $\Omega$ be a bounded connected open set in $\mathbb{R}^N$ ($N\ge 2$), with smooth enough boundary $\partial \Omega$. Let us introduce a set of periodically distributed holes. As a result, we obtain an open set $\Omega_\varepsilon$, where $\varepsilon$ represents a small parameter related to the characteristic size of the holes.

Let $Y=[0,1]^N$ be the representative cell in $\mathbb{R}^N$ and $F$ an open subset of $Y$ with smooth enough boundary $\partial F$, such that $\bar F\subset Y$. We denote $Y^*=Y\setminus \bar F$. For $k\in \mathbb{Z}^N$ and $\varepsilon\in (0,1]$, each cell $Y_{k,\varepsilon}=\varepsilon\,k+\varepsilon\,Y$ is similar to the unit cell $Y$ rescaled to size $\varepsilon$ and $F_{k,\varepsilon}=\varepsilon\,k+\varepsilon\,F$ is similar to $F$ rescaled to size $\varepsilon$. We denote $Y^*_{k,\varepsilon}=Y_{k,\varepsilon}\setminus \bar F_{k,\varepsilon}$. We denote by $F_\varepsilon$ the set of all the holes contained in $\Omega$, i.e. $F_\varepsilon=\displaystyle\cup_{k\in K}\{F_{k,\varepsilon}:\bar{F}_{k,\varepsilon}\subset \Omega \},$
where $K:=\{k\in \mathbb{Z}^N:Y_{k,\varepsilon}\cap\Omega\ne \emptyset \}$.

Let $\Omega_\varepsilon=\Omega\backslash \bar F_\varepsilon$. By this construction, $\Omega_\varepsilon$ is a periodically perforated domain with holes of the same size as the period. Remark that the holes do not intersect the boundary $\partial \Omega$. Let $\partial F_\varepsilon=\displaystyle\cup_{k\in K}\{\partial F_{k,\varepsilon}:\bar{F}_{k,\varepsilon}\subset \Omega \}.$ So
$$\partial \Omega_\varepsilon=\partial \Omega \cup \partial F_\varepsilon.$$
{\bf Position of the problem.} The prototype of the parabolic initial-boundary value problems that we consider in this article is
\begin{equation}
\left\{
\begin{array}
[c]{r@{\;}c@{\;}ll}%
\displaystyle\frac{\partial u_\varepsilon}{\partial t}-\Delta\,u_\varepsilon+\kappa u_\varepsilon &=&0
\quad & \text{\ in }\;\Omega_\varepsilon\times(0,T) ,\\
\nabla u_\varepsilon \cdot \nu+\varepsilon\,\displaystyle\frac{\partial u_\varepsilon}{\partial
t}&=&\varepsilon\,\delta \Delta_{\Gamma}u_\varepsilon-\varepsilon\,g(u_\varepsilon) & \text{\ on }%
\;\partial F_\varepsilon\times( 0,T),\\
u_\varepsilon&=& 0, & \text{\ on }%
\;\partial \Omega\times( 0,T),\\
u_\varepsilon(x,0) & = & u_\varepsilon^{0}(x), & \text{\ for }\;x\in\Omega_\varepsilon,\\
u_\varepsilon(x,0) & = & \psi_\varepsilon^{0}(x), & \text{\ for
}\;x\in\partial F_\varepsilon,
\end{array}
\right. \label{PDE}%
\end{equation}
where $u_\varepsilon=u_\varepsilon(x,t)$, $x\in \Omega_\varepsilon$, $t\in (0,T)$ and $T>0$. The first equation states the law of standard diffusion in $\Omega_\varepsilon$, $\Delta=\Delta_x$ denotes the Laplacian operator with respect to the space variable and $\kappa>0$ is a given constant. 
The boundary equation (\ref{PDE})$_{2}$ is multiplied by $\varepsilon$ to compensate the growth of the surface by shrinking $\varepsilon$, where the value of $u_\varepsilon$ is assumed to be the trace of the function $u_\varepsilon$ defined for $x\in \Omega_\varepsilon$, $\Delta_{\Gamma}$ denotes the Laplace-Beltrami operator on $\partial F_\varepsilon$, $\nu$ denotes the outward normal to $\partial F_\varepsilon$, and $\delta>0$ is a given constant. The term $\nabla u_\varepsilon \cdot \nu$ represents the interaction domain-boundary, while $\delta\Delta_{\Gamma}$ stands for a boundary diffusion. We assume that the function $g\in{C}\left(
\mathbb{R}\right) $ is given, and satisfies that there exist constants $q\ge 2$, $\alpha_1>0$, $\alpha_2>0,$
$\beta>0$, and $l>0$, such that
\begin{equation}
\alpha_{1}\left\vert s\right\vert ^{q}-\beta\leq g(s)s\leq\alpha
_{2}\left\vert s\right\vert ^{q}+\beta,\quad\text{for all
$s\in\mathbb{R}$,} \label{hip_2}%
\end{equation}
\begin{equation}
\left( g(s)-g(r)\right) \left( s-r\right) \geq-l\left(
s-r\right) ^{2},\quad\text{for all
$s,r\in\mathbb{R}$.} \label{hip_4}%
\end{equation}
Finally, we also assume that
\begin{equation}\label{hyp 0}
u_\varepsilon^{0}\in L^2\left( \Omega\right),\quad
\psi_\varepsilon^{0}\in L^{2}\left( \partial F_\varepsilon\right),
\end{equation}
are given, and we suppose that 
\begin{equation}\label{Initial_condition}
|u_\varepsilon^0|^2_{\Omega_\varepsilon}+\varepsilon|\psi_\varepsilon^0|^2_{\partial F_\varepsilon}\leq C,
\end{equation}
where $C$ is a positive constant, and we denote by $|\cdot|_{\Omega_\varepsilon}$ and $|\cdot|_{\partial F_\varepsilon}$ the norm in $L^2(\Omega_\varepsilon)$ and $L^2(\partial F_\varepsilon)$, respectively.

Depending of $\delta$, two classes of boundary conditions are modeled by (\ref{PDE}). For $\delta>0$, we have boundary conditions of {\it reactive-diffusive} type, and for $\delta=0$ the boundary conditions are {\it purely reactive}. In \cite{Anguiano}, we consider the homogenization of the problem (\ref{PDE}) with $\delta=0$ and we obtain rigorously a nonlinear parabolic problem with zero Dirichlet boundary condition and with extra-terms coming from the influence of the dynamical boundary conditions as the homogenized model. Though the results of the present article are similar to those of \cite{Anguiano}, the generalization of their proof is not trivial. Some new technical results are required in order to carry out the machinery of \cite{Anguiano}. Due to the presence of Laplace-Beltrami operator in the boundary condition, the variational formulation of the reaction-diffusion equation is different that in \cite{Anguiano}. We have to work in the space 
\begin{eqnarray}\label{space_W_delta}
W_\delta=\left\{ \left( v,\gamma_{0}(v)\right) \in H^1(\Omega_\varepsilon)\times H^1(\partial F_\varepsilon)\right\},\quad \delta>0,
\end{eqnarray}
where $\gamma_{0}$ denotes the trace operator $v\mapsto v|_{\partial\Omega_\varepsilon}$, and where we define by $H^1(\partial F_\varepsilon)$ the completion of $C^1(\partial F_\varepsilon)$ with respect to the induced norm by the inner product
$$((\phi, \psi))_{\partial F_\varepsilon}:=\int_{\partial F_\varepsilon}\phi\psi \,d\sigma+\delta \int_{\partial F_\varepsilon}\nabla_\Gamma \phi \cdot \nabla_\Gamma \psi d\sigma,\quad \forall \phi,\psi\in C^1(\partial F_\varepsilon),$$
where $\nabla_\Gamma$ denotes the tangential gradient on $\partial F_\varepsilon$ and $d\sigma$ denotes the natural volume element on $\partial F_\varepsilon$. The estimates of \cite{Anguiano} did not allow to cover this case and new estimates are needed to deal with problem (\ref{PDE}). In order to prove estimates in $H^2$-norm, we have to combine estimates for general elliptic boundary value problems with interpolation properties of Sobolev spaces (see Lemma \ref{estimates_NEW}). On the other hand, in order to pass to the limit, as $\varepsilon \to 0$, for the term which involves the tangential gradient $\nabla_{\Gamma}$ we make use of a convergence result based on a technique introduced by Vanninathan \cite{Vanni} for the Steklov problem which transforms surface integrals into volume integrals. This convergence result can be used taking into account the estimates in $H^2$-norm. Several technical results are merely quoted, and we refer \cite{Anguiano} for their proof.
We present here a new result concerning the local problem, which involves the orthogonal projection (denoted by $P_\Gamma$), the tangential gradient (denoted by $\nabla_\Gamma$) and the tangential divergence (denoted by ${\rm div}_\Gamma$) on the boundary of the unit cell. More precisely, using the so-called energy method introduced by Tartar \cite{Tartar} and considered by many authors (see, for instance, Cioranescu and Donato \cite{Ciora2}), we prove the following:

\begin{theorem}[Main Theorem]\label{Main}
Under the assumptions (\ref{hip_2})--(\ref{hip_4}) and (\ref{Initial_condition}), assume that $g\in\mathcal{C}^{1}( \mathbb{R})$, the exponent $q$ satisfies that 
\begin{equation}\label{assumption_q}
2\leq q<+\infty \  \text{ if } \ N=2 \quad \text{ and } \quad 2\leq q\leq {2N-2 \over N-2} \ \text{ if } \ N>2,
\end{equation}
and $(u_\varepsilon^0,\psi_\varepsilon^0)\in W_\delta\cap \left(L^{q}\left(
\Omega_\varepsilon\right) \times L^{q}\left( \partial F_\varepsilon\right)\right)$. Let $(u_\varepsilon, \psi_\varepsilon)$ be the unique solution of the problem (\ref{PDE}), where $\psi_\varepsilon(t)=\gamma_0(u_\varepsilon(t))$ a.e. $t\in (0,T]$. Then, as $\varepsilon\to 0$, we have
$$\tilde u_\varepsilon(t) \to u(t) \quad \text{strongly in } L^2(\Omega),\quad \forall t\in[0,T],$$ 
where $\tilde \cdot$ denotes the extension to $\Omega\times (0,T)$ and $u$ is the unique solution of the following problem
 \begin{equation}\label{limit_problem}
\left\{
\begin{array}{l}
\displaystyle \left({|Y^*|\over |Y|}+{|\partial F| \over |Y|} \right)\displaystyle\frac{\partial u}{\partial t}-{\rm div}\left(Q\nabla u \right)+ {|Y^*|\over |Y|}\kappa u+{|\partial F| \over |Y|}g(u)=0,
  \text{\ in }\;\Omega\times(0,T) ,\\[2ex]
u(x,0)  =  u_{0}(x),  \text{\ for }\;x\in\Omega,\\[2ex]
u= 0,  \text{\ on }
\;\partial \Omega\times( 0,T).
\end{array}
\right.
\end{equation}
The homogenized matrix $Q=((q_{i,j}))$, $1\leq i,j\leq N$, which is symmetric and positive-definite, is given by 
\begin{equation}\label{matrix}
q_{i,j}={1\over |Y|}\left(\int_{Y^*}\left(e_i+\nabla_y w_i \right)\cdot \left(e_j+\nabla_y w_j \right)dy+\delta\int_{\partial F}\left(P_\Gamma e_i+\nabla_\Gamma w_i \right)\cdot \left( P_\Gamma e_j+\nabla_\Gamma w_j\right)d\sigma(y)\right),
\end{equation}
where $w_i\in \mathbb{H}_{{\rm per}}/ \mathbb{R}$, $1\leq i\leq N$, is the unique solution of the cell problem 
\begin{equation}\label{system_eta}
\left\{
\begin{array}{l}
\displaystyle -{\rm div}_y\left(e_i+\nabla_y w_i \right)=0,   \text{\ in }Y^*,\\[2ex]
(e_i+\nabla_yw_i)\cdot \nu =\delta\,{\rm div}_{\Gamma}\left(P_{\Gamma}e_i+\nabla_\Gamma w_i \right),  \text{\ on }\partial F,\\[2ex]
w_i  \text{\ is }Y-\text{periodic}.
\end{array}
\right.
\end{equation}
Here, $e_i$ is the $i$ element of the canonical basis in $\mathbb{R}^N$ and $\mathbb{H}_{{\rm per}}$ is the space of functions from $\mathbb{H}:=\{v\in H^1(Y^*): v|_{\partial F}\in H^1(\partial F)\}$ which are $Y$-periodic.
\end{theorem}
\begin{remark}
Note that in the case $\delta=0$ (i.e., in the absence of a surface diffusion coefficient), the homogenized equation (\ref{limit_problem}) is exactly the equation obtained in \cite{Anguiano}.
\end{remark}
\begin{remark}
An example of a function $g\in\mathcal{C}^{1}( \mathbb{R})$, satisfying (\ref{hip_2})-(\ref{hip_4}), is a odd degree polynomial, 
$$g(s)=\sum_{j=0}^{2k+1}c_j\,s^j,$$
where $c_{2k+1}>0$. The typical case is $g(s)=s^3-s$.
\end{remark}
The homogenization of problems which involve the Laplace-Beltrami operator has been considered in recent articles. 

In particular, in the context of periodic homogenization based on the periodic unfolding method, in \cite{Graf_Peter} Graf and Peter extend the existing convergence results for the boundary periodic unfolding operator to gradients defined on manifolds. These results are then used to homogenize a system of five coupled reaction-diffusion equations, three of which include diffusion described by the Laplace-Beltrami operator and four of which consider a particular nonlinearity. 

In \cite{Amar_Gianni}, Amar and Gianni state a new property of the unfolding operator regarding the unfolded tangential gradient. This property is used to homogenize a differential system of linear equations in two disjoint conductive phases with a linear dynamical boundary condition which involves the Laplace-Beltrami operator in the separating interface. An error estimate for this model, under extra regularity assumptions on the data, can be found in Amar and Gianni \cite{Amar_Gianni2}. 

More recently, in \cite{Gahn}, Gahn derives some general two-scale compactness results for coupled bulk-surface problems and applies these results to an elliptic problem with a non-dynamical boundary condition, which involves the Laplace-Beltrami operator, in a multi-component domain. 

However, to our knowledge, there does not seem to be in the literature any study on the homogenization of parabolic models associated with nonlinear dynamical boundary conditions, which involves the Laplace-Beltrami operator, in a periodically perforated domain, as we consider in this article.

The article is organized as follows. In Section \ref{S2}, we introduce suitable functions spaces for our considerations. Especially, we consider some fundamentals from differential geometry as the tangential gradient and the tangential divergence. To prove the main result, in Section \ref{S3} we prove the existence and uniqueness of solution of (\ref{PDE}), {\it a priori} estimates are established in Section \ref{S4} and some compactness results are proved in Section \ref{S5}. Finally, the proof of Theorem \ref{Main} is established in Section \ref{S6}.

\section{Functional setting}\label{S2}
{\bf Notation.} We denote by $(\cdot,\cdot) _{\Omega_\varepsilon}$ (respectively, $(
\cdot,\cdot)_{\partial F_\varepsilon}$) the inner product in
$L^{2}(\Omega_\varepsilon)$ (respectively, in $L^{2}(\partial F_\varepsilon)$),
and by $\left\vert \cdot\right\vert _{\Omega_\varepsilon}$
(respectively, $\left\vert \cdot\right\vert
_{\partial F_\varepsilon}$) the associated norm. We also denote by $(\cdot,\cdot) _{\Omega_\varepsilon}$ the inner product in $(L^2(\Omega_\varepsilon))^N$.

If $r\ne2$, we will also denote by
$(\cdot,\cdot) _{\Omega_\varepsilon}$ (respectively, $(
\cdot,\cdot)_{\partial F_\varepsilon}$) the duality product between
$L^{r'}(\Omega_\varepsilon)$ and $L^{r}(\Omega_\varepsilon)$ (respectively, the duality
product between $L^{r'}(\partial F_\varepsilon)$ and
$L^{r}(\partial F_\varepsilon)$). We will denote by
$|\cdot|_{r,\Omega_\varepsilon}$ (respectively $|\cdot|_{r,\partial F_\varepsilon}$)
the norm in $L^r(\Omega_\varepsilon)$ (respectively in $L^r(\partial F_\varepsilon)$).

We denote by $(\cdot,\cdot) _{\Omega}$ the inner product in
$L^{2}(\Omega)$,
and by $\left\vert \cdot\right\vert _{\Omega}$ the associated norm. If $r\ne2$, we will also denote by
$(\cdot,\cdot) _{\Omega}$ the duality product between
$L^{r'}(\Omega)$ and $L^{r}(\Omega)$. We will denote by
$|\cdot|_{r,\Omega}$
the norm in $L^r(\Omega)$.

By $\left\Vert \cdot\right\Vert _{\Omega_\varepsilon}$ we denote the norm in
$H^{1}\left(\Omega_\varepsilon\right)$, which is associated to the inner
product $$((u,v))_{\Omega_\varepsilon}:=\left(u,v\right)
_{\Omega_\varepsilon}
+(\nabla u,\nabla v)_{\Omega_\varepsilon},\quad \forall u,v\in H^1(\Omega_\varepsilon),$$ and by $||\cdot||_{\Omega_\varepsilon,T}$ we denote the norm in $L^2(0,T;H^1(\Omega_\varepsilon))$. By $\left\Vert \cdot\right\Vert _{\Omega}$ we denote the norm in
$H^{1}\left(\Omega\right)$, by $||\cdot||_{\Omega,T}$ we denote the norm in $L^2(0,T;H^1(\Omega))$ and, if $r\ne2$, we denote by $|\cdot|_{r,\Omega,T}$ the norm in $L^r(0,T;L^r(\Omega))$.

We denote by $\gamma_{0}$ the trace operator $u\mapsto
u|_{\partial\Omega_\varepsilon}$, which belongs to
$\mathcal{L}(H^1(\Omega_\varepsilon), H^{1/2}(\partial\Omega_\varepsilon))$.

We introduce, for any $s> 1$, the space $H^s(\Omega_\varepsilon)$, which is naturally embedded in $H^1(\Omega_\varepsilon)$, and it is a Hilbert space equipped with the norm inherited, which we denote by $||\cdot ||_{H^s(\Omega_\varepsilon)}$.

Moreover, we denote by $H^r_{\partial \Omega}(\Omega_\varepsilon)$ and $H^r_{\partial \Omega}(\partial \Omega_\varepsilon)$, for $r\ge 0$, the standard Sobolev spaces which are closed subspaces of $H^r(\Omega_\varepsilon)$ and $H^r(\partial \Omega_\varepsilon)$, respectively, and the subscript $\partial \Omega$ means that, respectively, traces or functions in $\partial \Omega_\varepsilon$, vanish on this part of the boundary of $\Omega_\varepsilon$, i.e.
$$H^r_{\partial \Omega}(\Omega_\varepsilon)=\{v\in H^r(\Omega_\varepsilon):\gamma_0(v)=0  \text{ on } \partial \Omega \},$$
and
$$H^r_{\partial \Omega}(\partial \Omega_\varepsilon)=\{v\in H^r(\partial \Omega_\varepsilon):v=0  \text{ on } \partial \Omega \}.$$

Analogously, for $r\ge 2$, we denote
$$L^r_{\partial \Omega}(\partial \Omega_\varepsilon):=\{v\in L^r(\partial \Omega_\varepsilon):v=0  \text{ on } \partial \Omega \}.$$
Let us notice that, in fact, we can consider the given $\psi_\varepsilon^{0}$ as an element of $L^2_{\partial \Omega}(\partial \Omega_\varepsilon)$.

Let us consider the space
$$H_q:= L^{q}\left(
\Omega_\varepsilon\right) \times L_{\partial \Omega}^{q}\left( \partial\Omega_\varepsilon\right)
\text{,}\quad \forall q\ge 2,
$$with the natural inner product $ ((
v,\phi), ( w,\varphi))_{H_q}=(v,w)_{\Omega_\varepsilon}+
\varepsilon(\phi,\varphi)_{\partial F_\varepsilon},$ which in particular
induces the norm $|(\cdot,\cdot)|_{H_q}$ given by
$$|\left(
v,\phi\right)|^q_{H_q}=|v|_{q,\Omega_\varepsilon}^q+\varepsilon|\phi|^q_{q,\partial F_\varepsilon},\quad(v,\phi)\in
H_q.$$
For the sake of clarity, we shall omit to write explicitly the index $q$ if $q=2$, so we denote by $H$ the Hilbert space $$H:=L^{2}\left(
\Omega_\varepsilon\right) \times L_{\partial \Omega}^{2}\left( \partial\Omega_\varepsilon\right).$$
For functions $u\in H_{\partial \Omega}^1(\Omega_\varepsilon)$ which satisfy $\Delta u\in L_{\partial \Omega}^2(\Omega_\varepsilon)$, we have
$$\int_{\Omega_\varepsilon}\Delta u\,vdx=-\int_{\Omega_\varepsilon}\nabla u\cdot \nabla v dx+\int_{\partial F_\varepsilon}\nabla u \cdot \nu vd\sigma(x),\quad \forall v\in H_{\partial \Omega}^1(\Omega_\varepsilon).$$
{\bf Tangential gradient and Laplace-Beltrami operator.} We recall here, for the reader's convenience, some well-known facts on the tangential gradient $\nabla_\Gamma$ and the Laplace-Beltrami operator $\Delta_\Gamma$. We refer to Sokolowski and Zolesio \cite{Soko_Zole} for more details and proofs.

Let $S$ be a smooth surface with normal unit vector $\nu$. For every $v\in (L^2(S))^N$, we can define an element $P_{\Gamma}v\in (L^2(S))^N$ such that $P_{\Gamma}v\cdot \nu=0$ a.e. on $S$, where $P_\Gamma(y)$ for $y\in S$ is the orthogonal projection on the tangent space at $y\in S$, i.e., it holds that
$$P_\Gamma (y)v(y)=v(y)-\left(v(y)\cdot \nu(y)\right)\nu(y)\quad \text{for a.e. }y\in S.$$

Let $\phi \in C^1(S)$, there exist a tubular neighborhood $U$ of $S$ and an extension $\tilde \phi \in C^1(U)$ of $\phi$. We define the tangential gradient of $\phi$ on $S$ by
$$\nabla_\Gamma\phi:=P_\Gamma \nabla\tilde \phi=\nabla \tilde \phi-(\nabla \tilde \phi\cdot \nu)\nu \quad \text{on } S.$$
We emphasize that this definition is independent of the chosen extension of $\phi$.

Let $\Phi\in (C^1(S))^N$, then there exists an extension $\tilde \Phi \in (C^1(U))^N$ ($U$ as above a suitable neighborhood of $S$) and we define the tangential divergence of $\Phi$ on $S$ by
$${\rm div}_\Gamma \Phi:=\nabla_\Gamma \cdot \Phi:=\nabla\cdot \tilde \Phi-D\tilde \Phi\nu\cdot \nu \quad \text{on }S,$$
where $D\tilde \Phi$ is the Jacobi-matrix of $\tilde \Phi$. 

Now, we consider the surface $\partial F_\varepsilon$. First, an equivalent definition of the Sobolev space $H^1(\partial F_\varepsilon)$ on $\partial F_\varepsilon$ is given. We introduce the inner product
$$((\phi, \psi))_{\partial F_\varepsilon}:=(\phi,\psi)_{\partial F_\varepsilon}+\delta(\nabla_\Gamma \phi,\nabla_\Gamma \psi)_{\partial F_\varepsilon},\quad  \forall \phi,\psi\in C^1(\partial F_\varepsilon),\quad \delta\ge0,$$
and denote by $||\cdot ||_{\partial F_\varepsilon}$ the induced norm. The Sobolev space $H^1(\partial F_\varepsilon)$ is the closure of the space $C^1(\partial F_\varepsilon)$ with respect to the norm induced by the inner product. Therefore, the space $C^1(\partial F_\varepsilon)$ is dense by definition in the space $H^1(\partial F_\varepsilon)$. An equivalent definition of $H^1(\partial F_\varepsilon)$ can be given via local coordinates or distributional meaning, see, for instance, Strichartz \cite{Strichartz}. We denote by $||\cdot||_{\partial F_\varepsilon,T}$ the norm in $L^2(0,T;H^1(\partial F_\varepsilon))$.

By definition, for every $\phi \in H^1(\partial F_\varepsilon)$ there exists $\nabla_\Gamma \phi \in L^2(\partial F_\varepsilon)$ with $\nabla_\Gamma \phi \cdot \nu=0$ a.e. on $\partial F_\varepsilon$, the tangential gradient in the distributional sense. 

We introduce, for any $s> 1$, the space $H^s(\partial F_\varepsilon)$, which is naturally embedded in $H^1(\partial F_\varepsilon)$, equipped with the norm inherited, which we denote by $||\cdot ||_{H^s(\partial F_\varepsilon)}$.

For all $\psi \in H^1(\partial F_\varepsilon)$ and $v\in (C^1(\partial F_\varepsilon))^N$ such that $v\cdot \nu=0$ a.e. on $\partial F_\varepsilon$, we have the Stokes formula (see \cite[Proposition 2.58]{Soko_Zole})
\begin{equation}\label{Stokes_formula}
\int_{\partial F_\varepsilon}\nabla_\Gamma \psi\cdot v\,d\sigma=-\int_{\partial F_\varepsilon}\psi {\rm div}_\Gamma v \,d\sigma.
\end{equation}
Let $h\in H^2(\partial F_\varepsilon)$, then we have $\nabla_\Gamma h\in H^1(\partial F_\varepsilon)$ such that $\nabla_\Gamma h \cdot \nu=0$ a.e. on $\partial F_\varepsilon$. The Laplace-Beltrami operator $\Delta_\Gamma$ on $\partial F_\varepsilon$ is defined as follows
$$\Delta_\Gamma h={\rm div}_\Gamma\left(\nabla_\Gamma h \right)\quad \forall h\in H^2(\partial F_\varepsilon).$$
Hence $\Delta_\Gamma h\in L^2(\partial F_\varepsilon)$, and from (\ref{Stokes_formula}) it follows that the element $\Delta_\Gamma h\in L^2(\partial F_\varepsilon)$ is uniquely determined by the integral identity
\begin{equation}\label{Beltrami_formula}
\int_{\partial F_\varepsilon}\Delta_\Gamma h \psi\,d\sigma=-\int_{\partial F_\varepsilon}\nabla_\Gamma h\cdot \nabla \psi d\sigma\quad \forall \psi\in H^1(\partial F_\varepsilon).
\end{equation}
If $\psi\in H^1(\partial F_\varepsilon)$, then there exists (see \cite[Chapter 2, Section 2.20]{Soko_Zole}) an element $\vartheta\in H^{3/2}(\Omega_\varepsilon)$, the extension of $\psi$, and
\begin{equation}\label{function_boundary}
\vartheta|_{\partial F_\varepsilon}=\psi,\quad \text{furthermore}\quad \nabla \vartheta\cdot \nu=0\quad \text{on  }\partial F_\varepsilon.
\end{equation}
Therefore $\nabla \vartheta=\nabla_\Gamma \psi$ on $\partial F_\varepsilon$. It should be noted that on the right-hand side of (\ref{Beltrami_formula}) there is the scalar product of vector fields $\nabla_\Gamma h$ and $\nabla_\Gamma \psi$ tangent to $\partial F_\varepsilon$.

On the other hand, if $\psi$ is a smooth function defined in an open neighbourhood of $\partial F_\varepsilon$ in $\Omega$, then (see \cite[Chapter 2, Section 2.20]{Soko_Zole})
$$\nabla_\Gamma h\cdot \left(\nabla \psi|_{\partial F_\varepsilon} \right)=\nabla_\Gamma h\cdot \nabla_\Gamma \psi$$
because of
$$\left(\nabla \psi \cdot \nu\nu\right)\cdot \nabla_\Gamma h=0.$$
Hence, if $\psi$ is the restriction to $\partial F_\varepsilon$ of a given function $\psi$ defined in $\Omega$, then
\begin{equation}\label{test_function_boundary}
\int_{\partial F_\varepsilon}\Delta_\Gamma h \psi\,d\sigma=-\int_{\partial F_\varepsilon}\nabla_\Gamma h\cdot \nabla \psi d\sigma\quad \forall \psi\in H^2(\Omega).
\end{equation}
{\bf The space $W_{\delta}$.} We now introduce, as anticipated in the introduction, the space $W_{\delta}$ given in (\ref{space_W_delta}) (see \cite[Subsection 2.2]{Gal} for more details). Let $V^\delta_{\partial \Omega}$, $\delta\ge 0$, be the completion of $C^1(\overline \Omega_{\varepsilon})$ in the norm
$$||u||^2_{V^{\delta}_{\partial \Omega}}:=\int_{\Omega_\varepsilon}\left(|u(x)|^2+ |\nabla u(x)|^2 \right)dx+\varepsilon\int_{\partial F_\varepsilon}\left(|u(x)|^2+\delta |\nabla_{\Gamma}u(x)|^2 \right)d\sigma(x).$$
Note that for any $f\in V^{\delta}_{\partial \Omega}$, we have $f\in H^1_{\partial \Omega}(\Omega_\varepsilon)$ so that $f_{\partial F_\varepsilon}$ makes sense in the trace sense. The space $V^{\delta}_{\partial \Omega}$ is topologically isomorphic to $H^1(\Omega_\varepsilon)\times H^1_{\partial \Omega}(\partial\Omega_\varepsilon)$ if $\delta>0$, and $V^{0}_{\partial \Omega}=H^1_{\partial \Omega}(\Omega_\varepsilon)$. 

For all $\delta\ge 0$, we define the linear space
$$W_\delta:=\left\{ \left( v,\gamma_{0}(v)\right):v \in V^{\delta}_{\partial \Omega} \right\}.$$ 

We emphasize that $W_\delta$ is not a product space as $V^{\delta}_{\partial \Omega}$. Clearly, $W_\delta\subset H$ densely since the trace operator acting on function $H^1(\Omega_\varepsilon)$ and into $H^{1/2}(\partial \Omega_\varepsilon)$ is bounded and onto, and $W_\delta$ is a Hilbert space with respect to the inner product inherited from $V^{\delta}_{\partial \Omega}$, $\delta \ge 0$. Thus, by definition we can identify
$$W_\delta=\left\{ \left( v,\gamma_{0}(v)\right) \in H^1(\Omega_\varepsilon)\times H^1_{\partial \Omega}(\partial\Omega_\varepsilon)\right\},\quad \text{ for each }\delta>0,$$
where
\begin{equation*}\label{def_norma_WDelta}
\left\Vert \left( v,\gamma_{0}(v)\right) \right\Vert^2 _{W_\delta}:=\left\Vert v\right\Vert^2 _{\Omega_\varepsilon}+\varepsilon \left\Vert \gamma_0(v)\right\Vert^2 _{\partial F_\varepsilon},
\end{equation*}
and
$$W_0=\left\{ \left( v,\gamma_{0}(v)\right) \in H^1(\Omega_\varepsilon)\times H^{1/2}_{\partial \Omega}(\partial\Omega_\varepsilon)\right\}.$$

\section{Existence and uniqueness of solution}\label{S3}
Along this paper, we shall denote by $C$ different constants which are independent of $\varepsilon$. We state in this section a result on the existence and uniqueness of solution of problem (\ref{PDE}). 
First, we observe that it is easy to see from (\ref{hip_2}) that there
exists a constant $C>0$
 such that%
\begin{equation}\label{hipo_consecuencia}
\left\vert g(s)\right\vert \leq C\left( 1+\left\vert
s\right\vert ^{q-1}\right),\quad\text{for all
$s\in\mathbb{R}$.}
\end{equation}

\begin{definition}\label{definition_weakSolution} A weak solution of (\ref{PDE}) is a pair of functions $(u_\varepsilon,\psi_\varepsilon)$, satisfying
\begin{equation}\label{weak0}
 u_\varepsilon\in
C([0,T];L^2(\Omega_\varepsilon)),\quad \psi_\varepsilon\in
C([0,T];L_{\partial \Omega}^2(\partial\Omega_\varepsilon)),\quad\hbox{
for all $T>0$,}
\end{equation}
\begin{equation}\label{weak1}
 u_\varepsilon\in L^2(0,T;H^1(\Omega_\varepsilon)),
 \quad\hbox{
for all $T>0$,}
\end{equation}
\begin{equation}\label{weak2}
\psi_\varepsilon\in L^2(0,T;H_{\partial \Omega}^{1}(\partial\Omega_\varepsilon))\cap
L^q(0,T;L_{\partial \Omega}^q(\partial\Omega_\varepsilon)),\quad\hbox{ for all $T>0$,}
\end{equation}
\begin{equation}\label{weak3}
 \gamma_0(u_\varepsilon(t))=\psi_\varepsilon(t),\quad\hbox{ a.e. $t\in (0,T],$}
 \end{equation}
 \begin{equation}\label{weak4}
\left\{
\begin{array}{l}
 \dfrac{d}{dt}(u_\varepsilon(t),v)_{\Omega_\varepsilon}+\varepsilon\,\dfrac{d}{dt}(
\psi_\varepsilon(t),\gamma_{0}(v))_{\partial F_\varepsilon}+(\nabla u_\varepsilon(t),\nabla v)_{
\Omega_\varepsilon}+\kappa(u_\varepsilon(t),v)_{\Omega_\varepsilon}\\[2ex]
+\varepsilon\,\delta (\nabla_{\Gamma} \psi_\varepsilon(t),\nabla_{\Gamma} \gamma_0(v))_{\partial F_\varepsilon} +\varepsilon\,(g(\psi_\varepsilon(t)),\gamma_{0}%
(v))_{\partial F_\varepsilon}
 =0\\[2ex]
 \hbox{in $\mathcal{D}'(0,T)$, for all $v\in H^1(\Omega_\varepsilon)$
 such that $\gamma_0(v)\in H_{\partial \Omega}^{1}(\partial\Omega_\varepsilon)\cap L_{\partial \Omega}^q(\partial \Omega_\varepsilon),$}
\end{array}
\right.
\end{equation}
\begin{equation}\label{weak5}
 u_\varepsilon(0)=u_\varepsilon^0,\quad and\quad \psi_\varepsilon(0)=\psi_\varepsilon^0.
\end{equation}
\end{definition}

We have the following result.

\begin{theorem}
\label{Existence_solution_PDE}Under the assumptions (\ref{hip_2})--(\ref{hip_4}) and (\ref{hyp 0}), there exists a unique solution
$(u_\varepsilon,\psi_\varepsilon)$
of the problem (\ref{PDE}). Moreover, this solution satisfies the
energy equality
\begin{eqnarray}
&&\frac{1}{2}\frac{d}{dt}\left(|(u_\varepsilon(t),\psi_\varepsilon(t))|^2_{H}\right)+|\nabla
u_\varepsilon(t)|^2_{\Omega_\varepsilon} +\kappa|u_\varepsilon(t)|^2_{\Omega_\varepsilon}+\varepsilon\,\delta |\nabla_{\Gamma}\psi_\varepsilon(t)|^2_{\partial F_\varepsilon} +\varepsilon\,(g(\psi_\varepsilon(t)),\psi_\varepsilon(t))_{\partial F_\varepsilon}=0,\label{energyequality}
\end{eqnarray}
a.e. $t\in(0,T).$
\end{theorem}
\begin{proof}
On the space $W_\delta$ we define a continuous symmetric linear
operator $A_\delta:W_\delta\rightarrow W_\delta^{\prime}$, given by
\begin{equation}\label{def_A1}
\langle A_\delta(( v,\gamma_{0}(v))) ,( w,\gamma _{0}(w))\rangle
=(\nabla v,\nabla w)_{\Omega_\varepsilon}+\kappa(v,w)_{\Omega_\varepsilon}+\varepsilon\,\delta (\nabla_{\Gamma} \gamma_0(v),\nabla_{\Gamma} \gamma_0(w))_{\partial F_\varepsilon}\text{, \
}\end{equation}
for all $(v,\gamma_0(v))$, $(w,\gamma_0(w))\in W_\delta$.

We observe that $A_\delta$ is coercive. In fact, for all $(v,\gamma_0(v))\in W_\delta$, we have
\begin{eqnarray*}\nonumber
\left\langle A_\delta \left(( v,\gamma_{0}(v))\right)
,\left( v,\gamma _{0}(v) \right) \right\rangle+|(v,\gamma_0(v))|^2_{H} &
\ge&
{\rm min}\,\{1,\kappa \}\left\Vert v\right\Vert _{\Omega_\varepsilon}^{2}+\varepsilon\,\delta |\nabla_{\Gamma}\gamma_0(v)|^2_{\partial F_\varepsilon}+|v|_{\Omega_\varepsilon}^2+\varepsilon|\gamma_0(v)|^2_{\partial F_\varepsilon}%
\\
& \ge &{\rm min}\,\{1,\kappa \}\left\Vert \left( v,\gamma_{0}(v)\right) \right\Vert^2 _{W_\delta}\text{.}\label{Coercitivity}
\end{eqnarray*}
Let us denote
$$ 
V_1=W_\delta,\quad A_1=A_\delta,\quad V_{2}=L^{2}\left( \Omega_\varepsilon\right) \times L_{\partial \Omega}^{q}\left(
\partial\Omega_\varepsilon\right),\quad A_{2}\left( v,\phi\right) =(0,\varepsilon\,g(\phi)).$$
From \eqref{hipo_consecuencia} one deduces that
$A_2:V_2\rightarrow V'_2$.

With this notation, and denoting $V=V_1\cap V_2,$ $p_1=2,$
$p_2=q,$ $\vec{u}_\varepsilon=(u_\varepsilon,\psi_\varepsilon)$, one has that
\eqref{weak0}--\eqref{weak5} is equivalent to
\begin{equation}\label{Weak1}
 \vec{u}_\varepsilon\in C([0,T];H),\quad \vec{u}_\varepsilon\in\bigcap_{i=1}^2 L^{p_i}(0,T;V_i),
 \quad\mbox{for all $T>0,$}
\end{equation}
\begin{equation}\label{Weak2}
(\vec{u}_\varepsilon)'(t)+\sum_{i=1}^2A_i(\vec{u}_\varepsilon(t))=0\quad\mbox{in
$\mathcal{D}'(0,T;V'),$}
\end{equation}
\begin{equation}\label{Weak3}
\vec{u}_\varepsilon(0)=(u_\varepsilon^0,\psi_\varepsilon^0).
\end{equation}

Applying a slight modification of
\cite[Chapter 2,Theorem 1.4]{Lions}, it is not difficult to see that problem
\eqref{Weak1}--\eqref{Weak3} has a unique solution. Moreover,
 $\vec{u}_\varepsilon$ satisfies the energy equality
$$\frac{1}{2}\frac{d}{dt}|\vec{u}_\varepsilon(t)|^2_H+\sum_{i=1}^2\langle A_i(\vec{u}_\varepsilon(t)),\vec{u}_\varepsilon(t)\rangle_i=
0\quad\mbox{a.e. $t\in(0,T),$}$$ where
 $\left\langle \cdot ,\cdot \right\rangle_i $ denotes the duality product
between $V_i^{\prime}$ and $V_i$. This last equality turns out to be just \eqref{energyequality}.
\end{proof}

\section{A priori estimates}\label{S4}
In this section we obtain some energy estimates for the solution of (\ref{PDE}). By (\ref{energyequality}) and taking into account (\ref{hip_2}), we have
\begin{eqnarray}\label{precont1}
&&\frac{d}{dt}\left(|(u_\varepsilon(t),\psi_\varepsilon(t))|^2_{H}\right)+2\min\left\{
1,\kappa\right\}\left\Vert
u_\varepsilon(t)\right\Vert _{\Omega_\varepsilon}^2+2\varepsilon \delta |\nabla_\Gamma \psi_\varepsilon(t)|^2_{\partial F_\varepsilon} +2\alpha_1\varepsilon\,|\psi_\varepsilon(t)|_{q,\partial F_\varepsilon}^q\leq
 2\beta\varepsilon\,|\partial F_\varepsilon|, 
\end{eqnarray}
where $|\partial F_\varepsilon|$ denotes the measure of $\partial F_\varepsilon$.

Observe that the number of holes is given by $$N(\varepsilon)={|\Omega| \over (2\varepsilon)^N}\left(1+o(1)\right),$$
then using the change of variable 
\begin{equation*}\label{dilatacion}
y={x \over \varepsilon},\quad d\sigma(y)=\varepsilon^{-(N-1)}d\sigma(x),
\end{equation*}
we can deduce
\begin{equation}\label{acotaF}
|\partial F_\varepsilon|=N(\varepsilon) |\partial F_{k,\varepsilon}|=N(\varepsilon)\varepsilon^{N-1}|\partial F|\leq {C \over \varepsilon}.
\end{equation}
Let us denote
\[
\mathcal{G}(s):=\int_{0}^{s}g(r)dr.
\]
Then, there exist positive constants $\widetilde{\alpha}_{1}$,
$\widetilde {\alpha}_{2},$ and $\widetilde{\beta}$ such that
\begin{equation}
\widetilde{\alpha}_{1}| s|^{q}-\widetilde{\beta} \leq\mathcal{G}(s)\leq\widetilde{\alpha}_{2}|s|^{q}+\widetilde{\beta}\quad\forall s\in\mathbb{R}. \label{hip_2_adicional}
\end{equation}

We observe that the linear term $ \Delta_{\Gamma}u_\varepsilon$ in the boundary condition is coercive, so that this term is of no real significance to the energy estimates and only enhances the regularity of the solution. 
\begin{lemma}\label{estimates3}
Under the assumptions (\ref{hip_2})--(\ref{hip_4}) and (\ref{Initial_condition}), assume that $g\in\mathcal{C}^{1}( \mathbb{R})$. Then, for any initial condition $(u_\varepsilon^0,\psi_\varepsilon^0)\in W_\delta\cap H_q$, there exists a constant $C$ independent of $\varepsilon$, such that the solution $(u_\varepsilon,\psi_\varepsilon)$ of the problem (\ref{PDE}) satisfies
\begin{equation}\label{acotacion8_proof}
\left\Vert
u_\varepsilon\right\Vert _{\Omega_\varepsilon,T}\leq C, \quad  \sup_{t\in [0,T]}\left\Vert u_\varepsilon(t)\right\Vert _{\Omega_\varepsilon}\leq C, \quad \left\Vert
u'_\varepsilon\right\Vert _{\Omega_\varepsilon,T}\leq C,\quad  \sqrt{\varepsilon}\left\Vert
\gamma_0(u'_\varepsilon)\right\Vert _{\partial F_\varepsilon,T}\leq C,
\end{equation}
\begin{equation}\label{new_acotacion}
\sqrt{\varepsilon}|\psi_\varepsilon(t)|_{\partial F_\varepsilon}\leq C,\quad |u'_\varepsilon(t)|_{\Omega_\varepsilon}\leq C,\quad\sqrt{\varepsilon}|\gamma_0(u'_\varepsilon(t))|_{\partial F_\varepsilon}\leq C,
\end{equation}
for all $t\in (0,T)$.
\end{lemma}
\begin{proof}
Taking into account (\ref{acotaF}) in (\ref{precont1}), in particular, we obtain
\begin{eqnarray}\label{precont1_new}
&&\frac{d}{dt}\left(|(u_\varepsilon(t),\psi_\varepsilon(t))|^2_{H}\right)+2\min\left\{
1,\kappa\right\}\left\Vert
u_\varepsilon(t)\right\Vert _{\Omega_\varepsilon}^2\leq
C.
\end{eqnarray}
Integrating between $0$ and $t$ and taking into account (\ref{Initial_condition}), we obtain the first estimate in (\ref{acotacion8_proof}) and the fist estimate in (\ref{new_acotacion}).

Now, if we want to take the inner product in (\ref{PDE}) with $u'_\varepsilon$, we need that $u'_\varepsilon \in L^2(0,T;H^1(\Omega_\varepsilon))$ with $\gamma_0(u'_\varepsilon)\in L^2(0,T;H_{\partial \Omega}^{1}(\partial\Omega_\varepsilon))\cap
L^q(0,T;L_{\partial \Omega}^q(\partial\Omega_\varepsilon))$. However, we do not have it for our weak solution. Therefore, we use the Galerkin method in order to prove, rigorously, new {\it a priori} estimates for $u_\varepsilon$.

Let us observe that the space $H^{1}(\Omega_\varepsilon)\times H_{\partial \Omega}^{1}(\partial\Omega_\varepsilon)$ is compactly imbedded in $H$, and
therefore, for the symmetric and coercive linear continuous operator $A_{\delta}:W_{\delta}\rightarrow W_{\delta}^{\prime}$, where $A_\delta$ is given by (\ref{def_A1}), there exists a non-decreasing sequence $0<\lambda_{1}\leq\lambda_{2}\leq\ldots$ of eigenvalues associated to the operator
$A_{\delta}$ with $\lim_{j\rightarrow\infty }\lambda_{j}=\infty,$ and there exists a Hilbert basis of $H$, $\{(w_{j},\gamma_{0}(w_{j})) :j\geq1\}$$\subset D(A_\delta)$, with $span\{(w_{j},\gamma_{0}(w_{j})):j\geq1\} $ densely embedded in $W_{\delta}$, such that
\[
A_{\delta}((w_{j},\gamma_{0}(w_{j})))=\lambda_{j}(w_{j},\gamma_{0}(w_{j}))\quad\forall j\geq1.
\]

Taking into account the above facts, we denote by $$(u_{\varepsilon,m}(t),\gamma_{0}(u_{\varepsilon,m}(t)))=(u_{\varepsilon,m}(t;0,u_\varepsilon^0,\psi_\varepsilon^0),\gamma_{0}(u_{\varepsilon,m}(t;0,u_\varepsilon^0,\psi_\varepsilon^0)))$$ the Galerkin approximation of the solution
$(u_\varepsilon(t;0,u_\varepsilon^0,\psi_\varepsilon^0),\gamma_{0}(u_\varepsilon(t;0,u_\varepsilon^0,\psi_\varepsilon^0)))$ to (\ref{PDE}) for each integer $m\geq1$, which is given by
\begin{equation}
(u_{\varepsilon,m}(t),\gamma_{0}(u_{\varepsilon,m}(t)))=\sum_{j=1}^{m}\delta_{\varepsilon mj}(t)(w_{j},\gamma_{0}(w_{j})),\label{Galerkin1}
\end{equation}
and is the solution of
\begin{eqnarray}
\nonumber &&\dfrac{d}{dt}((u_{\varepsilon,m}(t),\gamma_{0}(u_{\varepsilon,m}(t))),(w_{j},\gamma_{0}(w_{j})))_{H}+\left\langle A_{\delta}((u_{\varepsilon,m}(t),\gamma_{0}(u_{\varepsilon,m}(t)))),(w_{j},\gamma_{0}(w_{j}))\right\rangle\\
&&+\varepsilon(  g(\gamma_{0}(u_{\varepsilon,m}(t))),\gamma_{0}(w_{j}))_{\partial F_\varepsilon}
=0,\quad j=1,\ldots ,m,\label{7}
\end{eqnarray}
with initial data
\begin{equation}
\label{7'}
(u_{\varepsilon,m}(0),\gamma_{0}(u_{\varepsilon,m}(0)))=(u_{\varepsilon,m}^{0},\gamma_{0}(u_{\varepsilon,m}^{0})),
\end{equation}
where
\[
\delta_{\varepsilon mj}(t)=(u_{\varepsilon,m}(t),w_{j})_{\Omega_\varepsilon}+( \gamma_{0}(u_{\varepsilon,m}(t)),\gamma_{0}(w_{j}))_{\partial F_\varepsilon},
\]
and $(u_{\varepsilon,m}^{0},\gamma_0(u_{\varepsilon,m}^{0}))\in span\{(w_j,\gamma_0(w_j)): j=1,\ldots ,m\}$ converge (when
$m\to\infty$) to $(u_\varepsilon^0,\psi_\varepsilon^0)$ in a suitable sense which will be specified below.

Let $(u_\varepsilon^0,\psi_\varepsilon^0)\in W_{\delta}\cap H_q$. For all $m\geq1$, since $span\{(w_{j},\gamma_{0}(w_{j})):j\geq1\} $ is densely embedded in $W_{\delta}\cap H_q$, there exists $(u_{\varepsilon,m}^{0},\gamma_{0}(u_{\varepsilon,m}^{0}))\in span\{(w_{j},\gamma_{0}(w_{j})):1\leq
j\leq m\} $, such that the sequence $\{(u_{\varepsilon,m}^{0},\gamma_{0}(u_{\varepsilon,m}^{0}))\}$ converges to $(u_\varepsilon^0,\psi_\varepsilon^0)$ in $W_{\delta}$
and in $H_q$. Then, in particular we know that there exists a constant $C$ such that
\begin{equation}\label{Initial_condition2}
||(u_{\varepsilon,m}^{0},\gamma_{0}(u_{\varepsilon,m}^{0}))||_{W_\delta}\leq C,\quad | (u_{\varepsilon,m}^{0},\gamma_{0}(u_{\varepsilon,m}^{0})) |_{H_q}\leq C.
\end{equation}
For each integer $m\geq1$, we consider the sequence $\{(u_{\varepsilon,m}(t),\gamma_{0}(u_{\varepsilon,m}(t)))\}$ defined by
(\ref{Galerkin1})-(\ref{7'}) with these initial data.

Multiplying by the derivative $\delta'_{\varepsilon mj}$ in (\ref{7}), and summing from $j=1$ to $m$, we obtain
\begin{eqnarray}\label{equality_G}
&& |(u_{\varepsilon,m}^{\prime}(t),\gamma_{0}(u_{\varepsilon,m}^{\prime}(t)))|^2_{H}+\frac{1}{2}\frac{d}{dt}(\left\langle A_{\delta}((u_{\varepsilon,m}(t),\gamma_{0}(u_{\varepsilon,m}(t)))),(u_{\varepsilon,m}(t),\gamma_{0}(u_{\varepsilon,m}(t)))\right\rangle)\nonumber\\
&&  +\varepsilon(g(\gamma_{0}(u_{\varepsilon,m}(t))),\gamma_{0}(u_{\varepsilon,m}^{\prime}(t)))_{\partial F_\varepsilon}=0.
\end{eqnarray}
We observe that
$$
(g(\gamma_{0}(u_{\varepsilon,m}(t))),\gamma_{0}(u_{\varepsilon,m}^{\prime}(t)))_{\partial F_\varepsilon}
=\frac{d}{dt}\int_{\partial F_\varepsilon}\mathcal{G}(\gamma_{0}(u_{\varepsilon,m}(t))) d\sigma(x).
$$
Then, integrating (\ref{equality_G}) between $0$ and $t$, taking into account the definition of $A_\delta$ and (\ref{acotaF})-(\ref{hip_2_adicional}), we obtain
\begin{eqnarray*}
&& \int_{0}^{t}|(u_{\varepsilon,m}^{\prime}(s),\gamma_{0}(u_{\varepsilon,m}^{\prime}(s)))|^2_{H}ds+ |\nabla
u_{\varepsilon,m}(t)|^2_{\Omega_\varepsilon} +\kappa|u_{\varepsilon,m}(t)|^2_{\Omega_\varepsilon}+\varepsilon\,\delta |\nabla_{\Gamma}\gamma_0(u_{\varepsilon,m}(t))|^2_{\partial F_\varepsilon} \nonumber\\
&&
+2\widetilde{\alpha}_{1}\varepsilon|\gamma_{0}(u_{\varepsilon,m}(t))|_{\partial F_\varepsilon}^{q} 
\leq\max\{1,\kappa\}\|(u_{\varepsilon,m}^{0},\gamma_{0}(u_{\varepsilon,m}^{0}))\|_{W_{\delta}}^{2}
+2\widetilde {\alpha}_{2}| (u_{\varepsilon,m}^{0}, \gamma_{0}(u_{\varepsilon,m}^{0}))|_{H_q}^{q} +4\tilde\beta C,
\end{eqnarray*}
for all $t\in (0,T)$, and we can deduce
\begin{eqnarray}\label{last_estimate}
&& \int_{0}^{t}|(u_{\varepsilon,m}^{\prime}(s),\gamma_{0}(u_{\varepsilon,m}^{\prime}(s)))|^2_{H}ds+ {\rm min}\,\{1,\kappa,2\tilde \alpha_1 \}||(u_{\varepsilon,m}(t),\gamma_0(u_{\varepsilon,m}(t)))||^2_{W_\delta}\nonumber\\
&\leq&C\left(1+\|(u_{\varepsilon,m}^{0},\gamma_{0}(u_{\varepsilon,m}^{0}))\|_{W_{\delta}}^{2}+| (u_{\varepsilon,m}^{0}, \gamma_{0}(u_{\varepsilon,m}^{0}))|_{H_q}^{q} \right),
\end{eqnarray}
for all $t\in (0,T)$. 
Taking into account (\ref{Initial_condition2}) in (\ref{last_estimate}), we have proved that the sequence $\{(u_{\varepsilon,m},\gamma_0(u_{\varepsilon,m}))\}$ is bounded in $C([0,T];W_\delta),$  and $\{(u_{\varepsilon,m}',\gamma_0(u_{\varepsilon,m}'))\}$ is bounded in $L^2(0, T;H),$
for all $T>0$.

If we work with the truncated Galerkin equations (\ref{Galerkin1})-(\ref{7'}) instead of the full PDE, we note that the calculations of the proof of  (\ref{precont1_new}) can be following identically to show that $\{u_{\varepsilon,m}\}$ is bounded in $L^2(0,T;H^1(\Omega_\varepsilon)),$
for all $T>0$.

Moreover, taking into account the uniqueness of solution to (\ref{PDE}) and using Aubin-Lions compactness lemma (e.g., cf. Lions \cite{Lions}), it is not difficult to conclude that the sequence $\{u_{\varepsilon,m}\}$ converges weakly in $
L^{2}(0,T;H^1(\Omega_\varepsilon))$ to the solution
$u_\varepsilon$ to (\ref{PDE}). Since the inclusion $H^1(\Omega_\varepsilon)\subset L^2(\Omega_\varepsilon)$ is compact and $u_\varepsilon\in C([0,T];L^2(\Omega_\varepsilon))$, it follows using \cite[Lemma 11.2]{Robinson} that  the second estimate in (\ref{acotacion8_proof}) is proved.


On the other hand, we note that, under the condition (\ref{hip_4}), we have that
\begin{equation}
g^{\prime}(s)\geq-l\quad\forall s\in\mathbb{R}. \label{47}
\end{equation}
Observe that as we are assuming that $g\in\mathcal{C}^{1}(\mathbb{R})$, we can differentiate with respect to time in (\ref{7}), and then, multiplying by the derivative $\delta'_{\varepsilon mj}$ and summing from $j=1$ to $m$, we obtain
\begin{eqnarray*}
&&  \frac{1}{2}\frac{d}{dt}|(u_{\varepsilon,m}^{\prime}(t),\gamma_{0}(u_{\varepsilon,m}^{\prime}(t)))|_{H}^{2}+\left\langle A_{\delta}((u_{\varepsilon,m}^{\prime}(t),\gamma_{0}(u_{\varepsilon,m}^{\prime}(t)))),(u_{\varepsilon,m}^{\prime}(t),\gamma_{0}(u_{\varepsilon,m}^{\prime}(t)))\right\rangle \\
&&=-\varepsilon(g^{\prime}(\gamma_{0}(u_{\varepsilon,m}(t))) \gamma_{0}(u_{\varepsilon,m}^{\prime}(t)),\gamma_{0}(u_{\varepsilon,m}^{\prime}(t)))_{\partial F_\varepsilon}.
\end{eqnarray*}
Then, using the definition of $A_\delta$ and (\ref{47}), we have
\begin{eqnarray*}
\frac{d}{dt}|(u_{\varepsilon,m}^{\prime}(t),\gamma_{0}(u_{\varepsilon,m}^{\prime}(t)))|_{H}^{2}+2|\nabla
u'_{\varepsilon,m}(t)|^2_{\Omega_\varepsilon} +2\kappa|u'_{\varepsilon,m}(t)|^2_{\Omega_\varepsilon}+2\varepsilon\,\delta |\nabla_{\Gamma}\gamma_0(u'_{\varepsilon,m}(t))|^2_{\partial F_\varepsilon}\leq 2l\varepsilon|\gamma_{0} (u_{\varepsilon,m}^{\prime}(t))|_{\partial F_\varepsilon}^{2},
\end{eqnarray*}
and we can deduce
\begin{eqnarray*}
&&\frac{d}{dt}|(u_{\varepsilon,m}^{\prime}(t),\gamma_{0}(u_{\varepsilon,m}^{\prime}(t)))|_{H}^{2}+2{\rm min}\,\{1,\kappa \}||u'_{\varepsilon,m}(t)||^2_{\Omega_\varepsilon}+2\varepsilon \left(|\gamma_{0} (u_{\varepsilon,m}^{\prime}(t))|_{\partial F_\varepsilon}^{2}+\delta |\nabla_{\Gamma}\gamma_0(u'_{\varepsilon,m}(t))|^2_{\partial F_\varepsilon} \right)\\
&\leq& 2\varepsilon(l+1)|\gamma_{0} (u_{\varepsilon,m}^{\prime}(t))|_{\partial F_\varepsilon}^{2}.
\end{eqnarray*}
Then, we obtain
\begin{eqnarray}\label{estimate_derivate}
\frac{d}{dt}|(u_{\varepsilon,m}^{\prime}(t),\gamma_{0}(u_{\varepsilon,m}^{\prime}(t)))|_{H}^{2}+2{\rm min}\,\{1,\kappa \}||(u'_{\varepsilon,m}(t),\gamma_0(u'_{\varepsilon,m}(t)))||^2_{W_\delta}
\leq 2(l+1)|(u_{\varepsilon,m}^{\prime}(t),\gamma_{0} (u_{\varepsilon,m}^{\prime}(t)))|_{H}^{2}.
\end{eqnarray}
Integrating between $r$ and $t$, we obtain
\begin{eqnarray*}
&&\!\!\!\!\!|(u_{\varepsilon,m}^{\prime}(t),\gamma_{0}(u_{\varepsilon,m}^{\prime}(t)))|_{H}^{2}+2{\rm min}\,\{1,\kappa \}\int_r^t ||(u'_{\varepsilon,m}(s),\gamma_0(u'_{\varepsilon,m}(s)))||^2_{W_\delta}ds \nonumber\\
 & \leq&| (  u_{\varepsilon,m}^{\prime}(r),\gamma_{0}(u_{\varepsilon,m}^{\prime}(r)))  |_{H}^{2} \!+\!2(l+1)\int_{r}^{t}|(u_{\varepsilon,m}^{\prime}(s),\gamma_{0}(u_{\varepsilon,m}^{\prime}(s)))|_{H}^{2}ds,
\end{eqnarray*}
for all $0\leq r\leq t$. Now, integrating with respect to $r$ between $0$ and $t$,
\begin{eqnarray*}
&&\!\!\!\!\!\!\!\!\!\!\!\!\!t|(u_{\varepsilon,m}^{\prime}(t),\gamma_{0}(u_{\varepsilon,m}^{\prime}(t)))|_{H}^{2}\!+2{\rm min}\,\{1,\kappa \}\int_0^t ||(u'_{\varepsilon,m}(s),\gamma_0(u'_{\varepsilon,m}(s)))||^2_{W_\delta}ds  \\
&  \leq&(2l+3)\int_{0}^{t}|(u_{\varepsilon,m}^{\prime}(s),\gamma_{0}(u_{\varepsilon,m}^{\prime}(s)))|_{H}^{2}ds,\nonumber
\end{eqnarray*}
for all $t\in(0,T)$, which, jointly with (\ref{last_estimate}), yields that
\begin{eqnarray}\label{lower1}
&&\int_0^t ||(u'_{\varepsilon,m}(s),\gamma_0(u'_{\varepsilon,m}(s)))||^2_{W_\delta}ds  
\leq C\left(1+\|(u_{\varepsilon,m}^{0},\gamma_{0}(u_{\varepsilon,m}^{0}))\|_{W_{\delta}}^{2}+| (u_{\varepsilon,m}^{0}, \gamma_{0}(u_{\varepsilon,m}^{0}))|_{H_q}^{q} \right),
\end{eqnarray}
and using (\ref{Initial_condition2}) we have proved that the sequence $\{(u_{\varepsilon,m}',\gamma_0(u_{\varepsilon,m}'))\}$ is bounded in $L^2(0, T;W_\delta),$ for all $T>0$. Then, the sequence $\{( u'_{\varepsilon,m},\gamma_{0} (u'_{\varepsilon, m}))\}$ converges weakly in $
L^{2}(0,T;W_\delta)$ to $(u'_\varepsilon,\gamma_{0}(u'_\varepsilon))$, for all $T>0$, and using the lower-semicontinuity of the norm and (\ref{lower1}), we get
\begin{eqnarray*}
||u'_\varepsilon||^2_{\Omega_\varepsilon,T}+\varepsilon ||\gamma_0(u'_\varepsilon)||^2_{\partial F_\varepsilon,T}&\leq& \liminf_{m\to \infty}\left(||u'_{\varepsilon,m}||^2_{\Omega_\varepsilon,T}+\varepsilon||\gamma_0(u'_{\varepsilon,m})||^2_{\partial F_\varepsilon,T} \right)\\
&\leq& C\liminf_{m\to \infty}\left(1+\|(u_{\varepsilon,m}^{0},\gamma_{0}(u_{\varepsilon,m}^{0}))\|_{W_{\delta}}^{2}+| (u_{\varepsilon,m}^{0}, \gamma_{0}(u_{\varepsilon,m}^{0}))|_{H_q}^{q} \right) \\
&=&C\left(1+\|(u_\varepsilon^0,\psi_\varepsilon^0)\|_{W_{\delta}}^{2}+| (u_\varepsilon^0,\psi_\varepsilon^0)|_{H_q}^{q} \right),
\end{eqnarray*}
which, jointly with $(u_\varepsilon^0,\psi_\varepsilon^0)\in W_\delta\cap H_q$, implies the last two estimates in (\ref{acotacion8_proof}).

On the other hand, for any $\tau>0$ and $t>\tau$, integrating (\ref{estimate_derivate}), in particular, we have
\begin{eqnarray*}
|(u_{\varepsilon,m}^{\prime}(r),\gamma_{0}(u_{\varepsilon,m}^{\prime}(r)))|_{H}^{2} 
   \leq | (  u_{\varepsilon,m}^{\prime}(\theta),\gamma_{0}(u_{\varepsilon,m}^{\prime}(\theta)))  |_{H}^{2} \!+\!2(l+1)\int_{\tau/2}^{t}|(u_{\varepsilon,m}^{\prime}(s),\gamma_{0}(u_{\varepsilon,m}^{\prime}(s)))|_{H}^{2}ds,
\end{eqnarray*}
for all $\tau/2\leq \theta\leq r\leq t$. Now, integrating with respect to $\theta$ between $\tau/2$ and $r$,
\begin{eqnarray*}
(r-\tau/2)|(u_{\varepsilon,m}^{\prime}(r),\gamma_{0}(u_{\varepsilon,m}^{\prime}(r)))|_{H}^{2} 
   \leq\left(2(l+1)(t-\tau/2)+1 \right)\int_{\tau/2}^{t}|(u_{\varepsilon,m}^{\prime}(s),\gamma_{0}(u_{\varepsilon,m}^{\prime}(s)))|_{H}^{2}ds,
\end{eqnarray*}
for all $0<\tau/2\leq r\leq t<T$, and, in particular
\begin{eqnarray*}
|(u_{\varepsilon,m}^{\prime}(r),\gamma_{0}(u_{\varepsilon,m}^{\prime}(r)))|_{H}^{2} 
   \leq2\tau^{-1}\left(2(l+1)(T-\tau/2)+1 \right)\int_{0}^{t}|(u_{\varepsilon,m}^{\prime}(s),\gamma_{0}(u_{\varepsilon,m}^{\prime}(s)))|_{H}^{2}ds,
\end{eqnarray*}
for all $r\in [\tau,t]$, which, jointly with (\ref{last_estimate}), yields that
\begin{eqnarray}\label{lower1_final}
|(u_{\varepsilon,m}^{\prime}(r),\gamma_{0}(u_{\varepsilon,m}^{\prime}(r)))|_{H}^{2} 
   \leq C\left(1+\|(u_{\varepsilon,m}^{0},\gamma_{0}(u_{\varepsilon,m}^{0}))\|_{W_{\delta}}^{2}+| (u_{\varepsilon,m}^{0}, \gamma_{0}(u_{\varepsilon,m}^{0}))|_{H_q}^{q} \right),
\end{eqnarray}
for all $r\in (0,T)$. 
Using (\ref{Initial_condition2}) we have proved that the sequence $\{(u_{\varepsilon,m}',\gamma_0(u_{\varepsilon,m}'))\}$ is bounded in $C([0,T];H)$. Then, the sequence $\{( u'_{\varepsilon,m}(r),\gamma_{0} (u'_{\varepsilon, m}(r)))\}$ converges weakly in $H$ to $(u'_\varepsilon(r),\gamma_{0}(u'_\varepsilon(r)))$, for all $r\in [0,T]$, and using the lower-semicontinuity of the norm and (\ref{lower1_final}), we get
\begin{eqnarray*}
|u'_\varepsilon(r)|^2_{\Omega_\varepsilon}+\varepsilon |\gamma_0(u'_\varepsilon(r))|^2_{\partial F_\varepsilon}&\leq& \liminf_{m\to \infty}\left(|u'_{\varepsilon,m}(r)|^2_{\Omega_\varepsilon}+\varepsilon|\gamma_0(u'_{\varepsilon,m}(r))|^2_{\partial F_\varepsilon} \right)\\
&\leq& C\liminf_{m\to \infty}\left(1+\|(u_{\varepsilon,m}^{0},\gamma_{0}(u_{\varepsilon,m}^{0}))\|_{W_{\delta}}^{2}+| (u_{\varepsilon,m}^{0}, \gamma_{0}(u_{\varepsilon,m}^{0}))|_{H_q}^{q} \right) \\
&=&C\left(1+\|(u_\varepsilon^0,\psi_\varepsilon^0)\|_{W_{\delta}}^{2}+| (u_\varepsilon^0,\psi_\varepsilon^0)|_{H_q}^{q} \right),
\end{eqnarray*}
which, jointly with $(u_\varepsilon^0,\psi_\varepsilon^0)\in W_\delta\cap H_q$, implies the last two estimates in (\ref{new_acotacion}).
\end{proof}
In the following result, we enhance the regularity of the solution.
\begin{lemma}\label{estimates_NEW}
Assume the assumptions in Lemma \ref{estimates3}. Then, for any initial condition $(u_\varepsilon^0,\psi_\varepsilon^0)\in W_\delta\cap H_q$, there exists a constant $C$ independent of $\varepsilon$, such that the solution $u_\varepsilon$ of the problem (\ref{PDE}) satisfies
\begin{equation}\label{acotacion_H2}
\left\Vert
 u_\varepsilon(t)\right\Vert _{H^2(\Omega_\varepsilon)}\leq C,
\end{equation}
for all $t\in (0,T)$.
\end{lemma}
\begin{proof}
In order to obtain the estimates for the $H^2$-norm, we rewrite (for every fixed $t$) problem (\ref{PDE}) as a second-order nonlinear elliptic boundary value problem:
\begin{equation}
\left\{
\begin{array}
[c]{r@{\;}c@{\;}ll}%
-\Delta\,u_\varepsilon+\kappa u_\varepsilon&=&h_1(t):=-\displaystyle\frac{\partial u_\varepsilon}{\partial t}
\quad & \text{\ in }\;\Omega_\varepsilon ,\\
-\varepsilon\,\delta \Delta_{\Gamma}u_\varepsilon+\varepsilon\,\lambda u_\varepsilon+\nabla u_\varepsilon \cdot \nu+\varepsilon\,g(u_\varepsilon)&=& \varepsilon\,h_2(t):=-\varepsilon \displaystyle\frac{\partial u_\varepsilon}{\partial
t}+\varepsilon\lambda u_\varepsilon \quad & \text{\ on }%
\;\partial F_\varepsilon,\\
u_\varepsilon&=& 0 & \text{\ on }%
\;\partial \Omega,
\end{array}
\right. \label{PDE_elliptic}%
\end{equation}
where $\lambda$ is some positive constant.

We multiply the first equation of (\ref{PDE_elliptic}) scalarly in $L^2(\Omega_\varepsilon)$ by $u_\varepsilon$, we integrate by parts and using (\ref{hip_2}), we have
\begin{eqnarray}\label{estimate_elliptic1}
|\nabla u_\varepsilon|^2_{\Omega_\varepsilon}+\kappa |u_\varepsilon|^2_{\Omega_\varepsilon}+\varepsilon\, \delta |\nabla_\Gamma \gamma_0(u_\varepsilon)|^2_{\partial F_\varepsilon}+\varepsilon\, \lambda | \gamma_0(u_\varepsilon)|^2_{\partial F_\varepsilon}+\varepsilon\,\alpha_1 |u_\varepsilon|^q_{q,\partial F_\varepsilon}\leq (h_1,u_\varepsilon)_{\Omega_\varepsilon}+\varepsilon (h_2,\gamma_0(u_\varepsilon))_{\partial F_\varepsilon}+\varepsilon\,\beta |\partial F_\varepsilon|.
\end{eqnarray}
Using Young's inequality, we obtain
\begin{eqnarray*}
(h_1,u_\varepsilon)_{\Omega_\varepsilon}\leq |h_1|_{\Omega_\varepsilon}|u_\varepsilon|_{\Omega_\varepsilon}\leq {1\over 2\kappa}|h_1|^2_{\Omega_\varepsilon}+{\kappa \over 2}|u_\varepsilon|^2_{\Omega_\varepsilon},
\end{eqnarray*}
and
\begin{eqnarray*}
(h_2,\gamma_0(u_\varepsilon))_{\partial F_\varepsilon}\leq |h_2|_{\partial F_\varepsilon}|\gamma_0(u_\varepsilon)|_{\partial F_\varepsilon}\leq {1\over 2\lambda}|h_2|^2_{\partial F_\varepsilon}+{\lambda \over 2}|\gamma_0(u_\varepsilon)|^2_{\partial F_\varepsilon},
\end{eqnarray*}
and by (\ref{estimate_elliptic1}), using (\ref{acotaF}), we can deduce, in particular, that there exists a positive constant $C$ such that
\begin{eqnarray}\label{estimate_elliptic2}
||u_\varepsilon||_{\Omega_\varepsilon}\leq C\left(1+|h_1|_{\Omega_\varepsilon}+\sqrt{\varepsilon} |h_2|_{\partial F_\varepsilon}\right).
\end{eqnarray}
Using now the estimates for general elliptic boundary value problems (see \cite[Chaper 2, Remark 7.2]{Lions_Magenes}) to the first equation of (\ref{PDE_elliptic}) with $s=2$, $m=1$ and $j=0$, we have
\begin{eqnarray}\label{regularity_elliptic1}
||u_\varepsilon||_{H^2(\Omega_\varepsilon)}\leq C\left(|h_1|_{\Omega_\varepsilon}+||\varepsilon\gamma_0(u_\varepsilon)||_{H^{3/2}(\partial F_\varepsilon)} \right).
\end{eqnarray}
Analogously, applying this estimate to the second equation in (\ref{PDE_elliptic}) and taking into account (\ref{acotaF}) and (\ref{47}), we deduce
\begin{eqnarray}\label{regularity_elliptic2}
||\varepsilon\gamma_0(u_\varepsilon)||_{H^2(\partial F_\varepsilon)}\leq C\left(1+\varepsilon |h_2|_{\partial F_\varepsilon}+| \partial_{\nu}u_\varepsilon|_{\partial F_\varepsilon} \right),
\end{eqnarray}
where by $\partial_{\nu}u_\varepsilon$ we denote $\nabla u_\varepsilon \cdot \nu$. Taking into account (\ref{regularity_elliptic2}) in (\ref{regularity_elliptic1}), we can deduce
\begin{eqnarray}\label{regularity_elliptic3}
||u_\varepsilon||_{H^2(\Omega_\varepsilon)}\leq C\left(1+|h_1|_{\Omega_\varepsilon}+\varepsilon |h_2|_{\partial F_\varepsilon}+| \partial_{\nu}u_\varepsilon|_{\partial F_\varepsilon} \right).
\end{eqnarray}
By the Trace Theorem in $H^{7/4}(\Omega_\varepsilon)$ (see \cite[Chapter 1, Theorem 9.4]{Lions_Magenes}), we have 
\begin{eqnarray*}
|\partial_{\nu}u_\varepsilon|_{\partial F_\varepsilon}\leq C||u_\varepsilon||_{H^{7/4}(\Omega_\varepsilon)},
\end{eqnarray*}
and by interpolation inequality (see \cite[Chapter 1, Remark 9.1]{Lions_Magenes}) with $s_1=1$, $s_2=2$ and $\theta=3/4$, we can deduce
\begin{eqnarray*}
|\partial_{\nu}u_\varepsilon|_{\partial F_\varepsilon}\leq C||u_\varepsilon||^{1/4}_{\Omega_\varepsilon}||u_\varepsilon||^{3/4}_{H^2(\Omega_\varepsilon)}.
\end{eqnarray*}
By Young's inequality, with the conjugate exponents $4$ and $4/3$, we get
\begin{eqnarray}\label{estimate_normal}
|\partial_{\nu}u_\varepsilon|_{\partial F_\varepsilon}\leq C\,||u_\varepsilon||_{\Omega_\varepsilon}+c\,||u_\varepsilon||_{H^2(\Omega_\varepsilon)},
\end{eqnarray}
where the positive constant $c$ can be arbitrarily small. Then, taking into account (\ref{estimate_normal}) in (\ref{regularity_elliptic3}), we have
\begin{eqnarray*}
||u_\varepsilon||_{H^2(\Omega_\varepsilon)}\leq C\left(1+|h_1|_{\Omega_\varepsilon}+\varepsilon |h_2|_{\partial F_\varepsilon}+||u_\varepsilon||_{\Omega_\varepsilon} \right),
\end{eqnarray*}
and using (\ref{estimate_elliptic2}), we can deduce the following estimate for the $H^2$-norm 
\begin{eqnarray}\label{estimate_H^2}
||u_\varepsilon||_{H^2(\Omega_\varepsilon)}\leq C\left(1+|h_1|_{\Omega_\varepsilon}+\sqrt{\varepsilon} |h_2|_{\partial F_\varepsilon} \right).
\end{eqnarray}
According to the second estimate in (\ref{new_acotacion}), we have
\begin{eqnarray}\label{estimate_h1}
|h_1|_{\Omega_\varepsilon}\leq C,
\end{eqnarray}
and by the first and third estimates in (\ref{new_acotacion}), we can deduce 
\begin{eqnarray}\label{estimate_h2}
\sqrt{\varepsilon} |h_2|_{\partial F_\varepsilon}\leq C.
\end{eqnarray}
Finally, taking into account (\ref{estimate_h1})-(\ref{estimate_h2}) in (\ref{estimate_H^2}), we obtain (\ref{acotacion_H2}).
\end{proof}

{\bf The extension of $u_\varepsilon$ to the whole $\Omega\times (0,T)$}: since the solution $u_\varepsilon$ of the problem (\ref{PDE}) is defined only in $\Omega_\varepsilon\times (0,T)$, we need to extend it to the whole $\Omega\times (0,T)$ to be able to state the convergence result. In order to do that, we use the well-known extension result given by Cioranescu and Saint Jean Paulin \cite{Cioranescu}. Taking into account Lemma \ref{estimates3}, the following result is a direct consequence of results contained in \cite[Corollary 4.8]{Anguiano}.
\begin{corollary}\label{estimates_extension}
Assume the assumptions in Lemma \ref{estimates3}. Then, there exists an extension $\tilde u_\varepsilon$ of the solution $u_\varepsilon$ of the problem (\ref{PDE}) into $\Omega\times (0,T)$, such that
\begin{equation}\label{acotacion1_extension}
\left\Vert
\tilde u_\varepsilon\right\Vert _{\Omega,T}\leq C, \quad |\tilde u_\varepsilon|_{q,\Omega,T}\leq C,
\end{equation}
\begin{equation}\label{acotacion2_extension}
\sup_{t\in [0,T]}\left\Vert \tilde u_\varepsilon(t)\right\Vert _{\Omega}\leq C, 
\end{equation}
\begin{equation}\label{acotacion3_extension}
 |\tilde u'_\varepsilon|_{q,\Omega,T}\leq C, 
\end{equation}
where the constant $C$ does not depend on $\varepsilon$.
\end{corollary}

 \section{A compactness result}\label{S5}
 In this section, we obtain some compactness results about the behavior of the sequence $\tilde u_\varepsilon$ satisfying the {\it a priori} estimates given in Corollary \ref{estimates_extension}.
 
 By $\chi_{\Omega_\varepsilon}$ we denote the characteristic function of the domain $\Omega_\varepsilon$. Due to the periodicity of the domain $\Omega_\varepsilon$, from Theorem 2.6 in Cioranescu and Donato \cite{CioraDonato} one has, for $\varepsilon \to 0$, that 
\begin{equation}\label{convergence_chi}
\chi_{\Omega_\varepsilon}\stackrel{\tt
*}\rightharpoonup {|Y^*|\over |Y|} \quad \textrm{weakly-star in}\
L^\infty(\Omega),
\end{equation}
where the limit is the proportion of the material in the cell $Y$.

Let $\xi_\varepsilon$ be the gradient of $u_\varepsilon$ in $\Omega_\varepsilon \times (0,T)$ and let us denote by $\tilde \xi_\varepsilon$ its extension with zero to the whole of $\Omega \times (0,T)$, i.e.
\begin{equation}\label{definition_tildexi}
\tilde \xi_\varepsilon=\left\{
\begin{array}{l}
\xi_\varepsilon \quad \text{in }\Omega_\varepsilon \times (0,T),\\
 0 \quad \text{in }(\Omega\setminus \overline{\Omega_\varepsilon})\times (0,T).
 \end{array}\right.
\end{equation}

\begin{proposition}\label{Propo_convergence}
Under the assumptions in Lemma \ref{estimates3}, there exists a function $u\in L^2(0,T;H_0^1(\Omega))\cap$\linebreak$ L^q(0,T;L^q(\Omega))$ ($u$ will be the unique solution of the limit system (\ref{limit_problem})) and a function $\xi\in L^2(0,T;L^2(\Omega))$ such that, at least after extraction of a subsequence, we have the following convergences for all $T>0,$
\begin{eqnarray}
\label{continuity1} &\tilde{u}_\varepsilon(t)\rightharpoonup
u(t) &\textrm{weakly in}\ H_0^1(\Omega),\quad \forall t\in[0,T],
\\
\label{converge_initial_data}
&\tilde{u}_{\varepsilon}(t)\rightarrow
u(t)&\quad \text{strongly in }L^2(\Omega),\quad \forall t\in[0,T],
\\
\label{converge_gradiente}
&\tilde \xi_\varepsilon \rightharpoonup
\xi &\quad \textrm{weakly in} \quad L^2(0,T;L^2(\Omega)),
\\
\label{convergence_gradiente2}
&\tilde \xi_\varepsilon \rightharpoonup
\xi &\quad \textrm{weakly in} \quad L^2(\Omega),\quad \forall t\in[0,T],
\end{eqnarray}
where $\tilde \xi_\varepsilon$ is given by (\ref{definition_tildexi}).

Let $q$ be the exponent satisfying (\ref{assumption_q}). Let $\bar q>1$ given by
$$
\bar q \in (1,2) \  \ \text{ if }N=2,\quad \bar q={2N\over (N-2)(q-2)+N}\ \ \text{ if } N>2.
$$
Then, we have the following convergences for all $T>0,$
\begin{eqnarray}
\label{converge_g_ae}
&g(\tilde u_{\varepsilon}(t))\rightarrow g(u(t))& \quad \text{strongly in} \quad L^{\bar q}(\Omega),\quad \forall t\in[0,T],
\\
\label{converge_g2_ae}
&g(\tilde u_{\varepsilon}(t))\rightharpoonup g(u(t))& \quad \text{weakly in} \quad W_0^{1,\bar q}(\Omega),\quad \forall t\in[0,T].
\end{eqnarray}
\end{proposition}
\begin{proof}
By (\ref{acotacion1_extension}), we observe that the sequence
$\{\tilde u_\varepsilon\}$ is bounded in the spaces $L^2(0,T;H_0^1(\Omega))\cap
 L^q(0,T;L^q(\Omega))$, for all $T>0.$ Let us fix $T>0$. Then, there exists a subsequence $\{\tilde u_{\varepsilon'}\}\subset \{\tilde u_\varepsilon\}$ and function $u\in L^2(0,T;H_0^1(\Omega))\cap L^q(0,T;L^q(\Omega))$ such that
\begin{eqnarray}
\label{continuity1_Rellich} &\tilde{u}_{\varepsilon'}\rightharpoonup
u &\textrm{weakly in}\ L^2(0,T;H_0^1(\Omega)),
\\
\label{continuity2} &\tilde{u}_{\varepsilon'}\rightharpoonup
u &\textrm{weakly in}\ L^q(0,T;L^q(\Omega)).
\end{eqnarray}
By the estimate (\ref{acotacion2_extension}), for each $t\in [0,T]$, we have that $\{\tilde u_\varepsilon(t)\}$ is bounded in $H_0^1(\Omega)$, and since we have (\ref{continuity1_Rellich}), we can deduce (\ref{continuity1}). By (\ref{continuity1}) and Rellich-Kondrachov Theorem, we obtain (\ref{converge_initial_data}). 

From the first estimate in (\ref{acotacion8_proof}) and (\ref{definition_tildexi}), we have $|\tilde \xi_\varepsilon|_{\Omega,T}\leq C$, and hence, up a sequence, there exists $\xi\in L^2(0,T,L^2(\Omega))$ such that 
\begin{eqnarray}\label{convergence_xi_debil}
&\tilde \xi_{\varepsilon''} \rightharpoonup \xi & \textrm{weakly in }  L^2(0,T;L^2(\Omega)),
\end{eqnarray}
and we have (\ref{converge_gradiente}). In order to prove (\ref{convergence_gradiente2}), we observe that by the estimate (\ref{acotacion2_extension}), for each $t\in [0,T]$, we have that $\tilde \xi_\varepsilon$ is bounded in $L^2(\Omega)$, and since we have (\ref{convergence_xi_debil}), we can deduce (\ref{convergence_gradiente2}).

 By the arbitrariness of $T>0$, all the convergences are satisfied, as we wanted to prove.

Now, we analyze the convergences for the nonlinear term $g$. By Rellich-Kondrachov Theorem, we have the compact embedding $H_0^1(\Omega)\subset L^r(\Omega)$ for all $r\in[2,2^{\star})$, where 
\begin{equation*}
2^{\star}=\left\{
\begin{array}{l}
{2N\over N-2} \quad \text{if }N>2,\\
 +\infty \quad \text{if }N=2.
 \end{array}\right.
\end{equation*}
By the estimate (\ref{acotacion2_extension}), for each $t\in[0,T]$, we have that $\{\tilde u_\varepsilon(t)\}$ is bounded in $H_0^1(\Omega)$. Then, the compact embedding $H_0^1(\Omega)\subset L^r(\Omega)$ for all $r\in[2,2^{\star})$, implies that it is precompact in $L^r(\Omega)$ for all $r\in[2,2^{\star})$.

By the estimate (\ref{acotacion3_extension}), we see that the sequence
$\{\tilde u'_\varepsilon\}$ is bounded in $L^r(0,T;L^r(\Omega))$, for all $T>0$ and for all $r\in[2,2^{\star})$.
Then, we have that $\tilde u_\varepsilon(t):[0,T]\longrightarrow L^r(\Omega)$ is an equicontinuous family of functions.

Then, applying the Ascoli-Arzelà Theorem, we deduce that $\{\tilde u_\varepsilon(t)\}$ is a precompact sequence in $C([0,T];L^r(\Omega))$ for all $r\in[2,2^{\star})$. Hence, since we have (\ref{continuity2}) for all $q\ge 2$, we can deduce that
\begin{eqnarray}\label{converge_strongly_p}
 &\tilde{u}_{\varepsilon'}\rightarrow
u &\textrm{strongly in}\ C([0,T];L^r(\Omega)),
\end{eqnarray}
for all $r\in[2,2^{\star})$.

We separate the cases $N>2$ and $N=2$.

{\bf Case 1: $N>2$.} Since
$$
{2^{\star}\over \bar q}={(N-2)(q-2)+N\over N-2}=q-1+{2\over N-2}>q-1,
$$
there exists $r\in [2,2^{\star})$ such that ${r\over \bar q}\ge q-1$ and
\begin{equation*}
\left\vert g(s)\right\vert \leq C\left( 1+\left\vert
s\right\vert ^{q-1}\right)\leq C\left( 1+\left\vert
s\right\vert ^{{r\over \bar q}}\right).
\end{equation*}
Then, applying Theorem 2.4 in \cite{Conca} for $G(x,v)=g(v)$, $t=\bar q$ and $r\in [2,2^{\star})$ such that ${r\over \bar q}\ge q-1$, we have that the map $v\in L^r(\Omega)\mapsto g(v)\in L^{\bar q}(\Omega)$ is continuous in the strong topologies. Then, taking into account (\ref{converge_strongly_p}), we get (\ref{converge_g_ae}).

Finally, we prove (\ref{converge_g2_ae}). First, we observe that it is easy to see from (\ref{hip_2}) that there
exists a constant $C>0$ such that
$$
\left\vert g'(s)\right\vert \leq C\left( 1+\left\vert
s\right\vert ^{q-2}\right).
$$
Then, we get
\begin{eqnarray}\label{inequality_CONCA}
\int_{\Omega}\left|{\partial g \over \partial x_i}(\tilde u_\varepsilon(t))\right|^{\bar q}dx&\leq& C\int_{\Omega}\left(1+|\tilde u_\varepsilon(t)|^{(q-2)\bar q} \right)\left|{\partial \tilde u_\varepsilon(t) \over \partial x_i}\right|^{\bar q}dx\\
&\leq& C\left( 1+\left(\int_{\Omega}|\tilde u_\varepsilon(t)|^{(q-2)\bar q \gamma}dx \right)^{1/\gamma} \right)\left(\int_{\Omega}|\nabla \tilde u_\varepsilon(t)|^{\bar q \eta} dx\right)^{1/\eta},\nonumber
\end{eqnarray}
where we took $\gamma$ and $\eta$ such that $\bar q \eta=2$, $1/\gamma+1/\eta=1$ and $(q-2)\bar q \gamma=2^{\star}$. Note that from here we get 
$$\bar q={2N\over (N-2)(q-2)+N}.$$ 
Observe that $\bar q>1$. Indeed,
$$
q\leq {2N-2\over N-2}={N\over N-2}+1< {N\over N-2}+2\Rightarrow (N-2)(q-2)+N<2N \Rightarrow {2N\over (N-2)(q-2)+N}>1.
$$
Then, we have
\begin{eqnarray*}
\int_{\Omega}\left|{\partial g \over \partial x_i}(\tilde u_\varepsilon(t))\right|^{\bar q}dx&\leq& C\left( 1+|\tilde u_\varepsilon|_{2^{\star},\Omega}^{2^{\star}/\gamma} \right)|\nabla \tilde u_\varepsilon|_{\Omega}^{2/ \eta},
\end{eqnarray*}
and taking into account the continuous embedding $H_0^1(\Omega)\subset L^{2^{\star}}(\Omega)$ and (\ref{acotacion2_extension}), we get
\begin{equation}\label{acotacion_Lqprime}
|\nabla g(\tilde u_\varepsilon(t))|_{\bar q,\Omega}\leq C.
\end{equation}
Then, from (\ref{converge_g_ae}) and (\ref{acotacion_Lqprime}), we can deduce (\ref{converge_g2_ae}).

{\bf Case 2: $N=2$.} We consider $s\in [2q-2,+\infty)$ and
\begin{equation}\label{def_barq2}
\bar q={2s\over 2(q-2)+s}.
\end{equation}
Since
$$
{s\over \bar q}={2(q-2)+s\over 2}=q-1+{s-2\over 2}\ge q-1,
$$
we have
\begin{equation*}
\left\vert g(s)\right\vert \leq C\left( 1+\left\vert
s\right\vert ^{q-1}\right)\leq C\left( 1+\left\vert
s\right\vert ^{{s\over \bar q}}\right).
\end{equation*}
Then, applying Theorem 2.4 in \cite{Conca} for $G(x,v)=g(v)$, $t=\bar q$ and $r=s$, we have that the map $v\in L^s(\Omega)\mapsto g(v)\in L^{\bar q}(\Omega)$ is continuous in the strong topologies. Then, taking into account (\ref{converge_strongly_p}), we get (\ref{converge_g_ae}).

Finally, we prove (\ref{converge_g2_ae}). In (\ref{inequality_CONCA}) we took $\gamma$ and $\eta$ such that $\bar q \eta=2$, $1/\gamma+1/\eta=1$ and $(q-2)\bar q \gamma=s$. Note that from here we get $\bar q$ given by (\ref{def_barq2}).

Observe that $\bar q\in (1,2)$. Indeed, taking into account that $\displaystyle{1\over \bar q}={q-2\over s}+{1\over 2}$, we can deduce
$$
2q-2\leq s<+\infty\Rightarrow0<{1\over s}\leq {1\over 2q-2}\Rightarrow {1\over 2}<{1\over \bar q}\leq {q-2\over 2q-2}+{1\over 2}\Rightarrow {2(2q-2)\over 2(q-2)+2q-2}\leq \bar q<2,
$$
and using that ${2(2q-2)\over 2(q-2)+2q-2}>1$, we have that $\bar q\in (1,2)$.

Then, we have
\begin{eqnarray*}
\int_{\Omega}\left|{\partial g \over \partial x_i}(\tilde u_\varepsilon(t))\right|^{\bar q}dx&\leq& C\left( 1+|\tilde u_\varepsilon|_{s,\Omega}^{s/\gamma} \right)|\nabla \tilde u_\varepsilon|_{\Omega}^{2/ \eta},
\end{eqnarray*}
and taking into account the continuous embedding $H_0^1(\Omega)\subset L^{s}(\Omega)$ and (\ref{acotacion2_extension}), we get
\begin{equation}\label{acotacion_LqprimeN2}
|\nabla g(\tilde u_\varepsilon(t))|_{\bar q,\Omega}\leq C.
\end{equation}
Then, from (\ref{converge_g_ae}) and (\ref{acotacion_LqprimeN2}), we can deduce (\ref{converge_g2_ae}).

\end{proof}
Because we have the linear term $ \Delta_{\Gamma}u_\varepsilon$ in the boundary condition, in order to pass to the limit in the integral which involves this term, we need the following result.
\begin{proposition}\label{Propo_convergence_New}
Under the assumptions in Lemma \ref{estimates_NEW}, there exists a function $\xi\in L^2(0,T;H^1(\Omega))$ such that for all $T>0,$
\begin{eqnarray}
\label{convergence_gradiente_new}
&\tilde \xi_\varepsilon \rightharpoonup
\xi &\quad \textrm{weakly in} \quad H^1(\Omega),\quad \forall t\in[0,T],
\end{eqnarray}
where $\tilde \xi_\varepsilon$ is given by (\ref{definition_tildexi}).
\end{proposition}
\begin{proof}
From the estimate (\ref{acotacion_H2}) and (\ref{definition_tildexi}), we have $||\tilde \xi_\varepsilon||_{\Omega}\leq C$. Then, we see that the sequence $\{\tilde \xi_\varepsilon\}$ is bounded in $H^1(\Omega)$, and hence, up to a subsequence and by (\ref{convergence_gradiente2}), we can deduce (\ref{convergence_gradiente_new}).
\end{proof}

\section{Homogenized model: proof of the main Theorem}\label{S6}
In this section, we identify the homogenized model.

We multiply system (\ref{PDE}) by a test function $v\in \mathcal{D}(\Omega)$, integrating by parts and taking into account (\ref{test_function_boundary}) and (\ref{definition_tildexi}), we have
\begin{eqnarray*}
\dfrac{d}{dt}\left(\int_{\Omega}\chi_{\Omega_\varepsilon}\tilde u_\varepsilon(t)vdx\right)+\varepsilon\,\dfrac{d}{dt}\left(\int_{\partial F_\varepsilon}
\gamma_{0}(u_\varepsilon(t))vd\sigma(x)\right)+\int_{\Omega} \tilde \xi_\varepsilon \cdot \nabla vdx
+\kappa \int_{\Omega}\chi_{\Omega_\varepsilon}\tilde u_\varepsilon(t)v dx\\[2ex]
+\varepsilon\, \delta\int_{\partial F_\varepsilon} \nabla_\Gamma \gamma_0(u_\varepsilon(t))\cdot\nabla vd\sigma(x)
 +\varepsilon\, \int_{\partial F_\varepsilon}
g(\gamma_{0}(u_\varepsilon(t)))vd\sigma(x)
 =0, 
 \end{eqnarray*}
 in $\mathcal{D}'(0,T)$.
 
 We consider $\varphi\in C_c^1([0,T])$ such that $\varphi(T)=0$ and $\varphi(0)\ne 0$. Multiplying by $\varphi$ and integrating between $0$ and $T$, we have
 \begin{eqnarray}\label{system1}\nonumber
-\varphi(0)\left(\int_{\Omega}\chi_{\Omega_\varepsilon}\tilde u_\varepsilon(0)vdx\right)-\int_0^T\dfrac{d}{dt}\varphi(t)\left(\int_{\Omega}\chi_{\Omega_\varepsilon}\tilde u_\varepsilon(t)vdx\right)dt\\[2ex]\nonumber
-\varepsilon\varphi(0)\left(\int_{\partial F_\varepsilon}
\gamma_{0}(u_\varepsilon(0))vd\sigma(x)\right)
-\varepsilon\int_0^T \dfrac{d}{dt}\varphi(t)\left(\int_{\partial F_\varepsilon}
\gamma_{0}(u_\varepsilon(t))vd\sigma(x)\right)dt\\[2ex]
+\int_0^T\varphi(t)\int_{\Omega}\tilde \xi_\varepsilon\cdot \nabla vdxdt+\kappa \int_0^T\varphi(t)\int_{\Omega}\chi_{\Omega_\varepsilon}\tilde u_\varepsilon(t)v dxdt\\[2ex]\nonumber
+\varepsilon\,\delta\int_0^T\varphi(t)\int_{\partial F_\varepsilon} \nabla_\Gamma \gamma_0(u_\varepsilon(t))\cdot \nabla vd\sigma(x)dt
 +\varepsilon \int_0^T \varphi(t)\int_{\partial F_\varepsilon}
g(\gamma_{0}(u_\varepsilon(t)))vd\sigma(x)dt=0.
  \end{eqnarray}
For the sake of clarity, we split the proof in three parts. Firstly, we pass to the limit, as $\varepsilon \to 0$, in (\ref{system1}) in order to get the limit equation satisfied by $u$. Secondly we identify $\xi$ making use of the solutions of the cell-problems (\ref{system_eta}), and finally we prove that $u$ is uniquely determined.

{\bf Step 1.} In order to pass to the limit, as $\varepsilon \to 0$, we reason as in \cite[Theorem 6.1]{Anguiano} for all the terms except the term which involves the tangential gradient $\nabla_{\Gamma}$. Exactly, for the integrals on $\Omega$ we only require to use Proposition \ref{Propo_convergence} and the convergence (\ref{convergence_chi}) and for the integrals on the boundary of the holes we make use of a convergence result based on a technique introduced by Vanninathan \cite{Vanni} for the Steklov problem which transforms surface integrals into volume integrals, which was already used as a main tool to homogenize the non homogeneous Neumann problem for the elliptic case by Cioranescu and Donato \cite{Ciora2}. For the term which involves the tangential gradient, we also use this technique together with Proposition \ref{Propo_convergence_New}. 

By Definition 3.2 in Cioranescu and Donato \cite{Ciora2}, let us introduce, for any $h\in L^{s'}(\partial F)$, $1\leq s'\leq \infty$, the linear form $\mu_{h}^\varepsilon$ on $W_0^{1,s}(\Omega)$ defined by
\begin{equation*}
\langle \mu_{h}^\varepsilon,\varphi \rangle=\varepsilon \int_{\partial F_\varepsilon} h\left(x\over \varepsilon \right)\varphi(x) d\sigma(x),\quad \forall \varphi\in W_0^{1,s}(\Omega),
\end{equation*}
with $1/s+1/s'=1$. It is proved in Lemma 3.3 in Cioranescu and Donato \cite{Ciora2} that
\begin{equation}\label{convergence_mu}
\mu_{h}^\varepsilon \to \mu_{h}\quad \text{strongly in }(W_0^{1,s}(\Omega))',
\end{equation}
where $\langle \mu_{h},\varphi \rangle=\displaystyle\mu_{h}\int_{\Omega}\varphi(x) dx$, with
$$\mu_{h}={1\over |Y|}\int_{\partial F}h(y)d\sigma(y).$$
In the particular case in which $h\in L^{\infty}(\partial F)$ or even when $h$ is constant, we have
\begin{equation*}
\mu_{h}^\varepsilon \to \mu_{h}\quad \text{strongly in }W^{-1,\infty}(\Omega).
\end{equation*}
We denote by $\mu_1^\varepsilon$ the above introduced measure in the particular case in which $h=1$. Notice that in this case $\mu_{h}$ becomes $\mu_1=|\partial F|/|Y|$.

For the term which involves the tangential gradient, we proceed as follows. Taking into account (\ref{function_boundary}), there exists an element $\vartheta_\varepsilon\in H^{3/2}(\Omega_\varepsilon)$, the extension of $\gamma_0(u_\varepsilon(t))$, such that $\nabla \vartheta_\varepsilon=\nabla_\Gamma \gamma_0(u_\varepsilon(t))$ on $\partial F_\varepsilon$. Then, we can deduce
\begin{eqnarray*}
\varepsilon\int_{\partial F_\varepsilon} \nabla_\Gamma \gamma_0(u_\varepsilon(t))\cdot \nabla vd\sigma(x)=\varepsilon\int_{\partial F_\varepsilon} \nabla \vartheta_\varepsilon\cdot \nabla vd\sigma(x)=\langle \mu_{1}^\varepsilon,\tilde\xi_\varepsilon \cdot \nabla v \rangle,
\end{eqnarray*}
where $\tilde \xi_\varepsilon$ is given by (\ref{definition_tildexi}).
Note that using (\ref{convergence_mu}) with $s=2$ and taking into account (\ref{convergence_gradiente_new}), we can deduce, for $\varepsilon\to 0$,
\begin{eqnarray*}
\varepsilon\int_{\partial F_\varepsilon} \nabla_\Gamma \gamma_0(u_\varepsilon(t))\cdot \nabla vd\sigma(x)=\langle \mu_{1}^\varepsilon,\tilde\xi_\varepsilon \cdot\nabla v \rangle\to \mu_1\int_{\Omega}\xi\cdot \nabla v dx={|\partial F|\over |Y|}\int_{\Omega}\xi\cdot \nabla vdx,\quad \forall v\in \mathcal{D}(\Omega),
\end{eqnarray*}
which integrating in time and using Lebesgue's Dominated Convergence Theorem, gives
\begin{eqnarray}\label{limit_Beltrami}
\varepsilon \int_0^T\varphi(t)\int_{\partial F_\varepsilon} \nabla_\Gamma \gamma_0(u_\varepsilon(t))\cdot \nabla vd\sigma(x)dt\to {|\partial F|\over |Y|}\int_0^T\varphi(t)\left( \int_{\Omega}\xi\cdot \nabla vdx\right)dt.
\end{eqnarray}
Therefore, using the proof of the main Theorem in \cite{Anguiano} and (\ref{limit_Beltrami}), we pass to the limit, as $\varepsilon \to 0$, in (\ref{system1}), and we obtain
 \begin{eqnarray*}
-\varphi(0)\left({|Y^*|\over |Y|}+{|\partial F| \over |Y|} \right)\left(\int_{\Omega} u(0)vdx\right)-\left({|Y^*|\over |Y|}+{|\partial F| \over |Y|} \right)\int_0^T\dfrac{d}{dt}\varphi(t)\left(\int_{\Omega} u(t)vdx\right)dt\\[2ex]\nonumber
+\int_0^T\varphi(t)\int_{\Omega}\xi\cdot \nabla vdxdt+\kappa {|Y^*|\over |Y|}\int_0^T\varphi(t)\int_{\Omega} u(t)v dxdt\\[2ex]\nonumber
+\delta {|\partial F|\over |Y|}\int_0^T\varphi(t)\int_{\Omega}\xi\cdot \nabla vdxdt
 +{|\partial F| \over |Y|} \int_0^T \varphi(t)\int_{\Omega}
g(u(t))vdxdt=0.
\end{eqnarray*}
 Hence, $\xi$ verifies
 \begin{equation}\label{equation_xi}
 \left({|Y^*|\over |Y|}+{|\partial F| \over |Y|} \right)\displaystyle\frac{\partial u}{\partial t}-\left(1+\delta{|\partial F|\over |Y|}\right){\rm div}\xi+ {|Y^*|\over |Y|}\kappa u+{|\partial F| \over |Y|}g(u)= 0, \quad \text{in }\Omega\times (0,T).
 \end{equation}
  
{\bf Step 2.} It remains now to identify $\xi$. We shall make use of the solutions of the cell problems (\ref{system_eta}). For any fixed $i=1,...,N$, let us define
\begin{equation}\label{definition_Psi}
\Psi_{i\varepsilon}(x)=\varepsilon\left(w_i\left({x\over \varepsilon}\right)+y_i \right)\quad \forall x\in \Omega_\varepsilon,
\end{equation}
where $y=x/\varepsilon$.

By periodicity
\begin{equation*}\label{convergence_Psi}
\tilde \Psi_{i\varepsilon}\rightharpoonup
x_i \quad \textrm{weakly in} \quad H^1(\Omega),
\end{equation*}
where $\tilde \cdot$ denotes the extension to $\Omega$ given by Cioranescu and Saint Jean Paulin \cite{Cioranescu}. Then, by Rellich-Kondrachov Theorem, we can deduce
\begin{equation}\label{convergence_Psi_fuerte}
\tilde \Psi_{i\varepsilon}\rightarrow
x_i \quad \textrm{strongly in} \quad L^2(\Omega).
\end{equation}
Let $\nabla \Psi_{i\varepsilon}$ be the gradient of $\Psi_{i\varepsilon}$ in $\Omega_\varepsilon$. Denote by $\widetilde{\nabla \Psi_{i\varepsilon}}$ the extension by zero of $\nabla \Psi_{i\varepsilon}$ inside the holes. From (\ref{definition_Psi}), we have
$$\widetilde{\nabla \Psi_{i\varepsilon}}=\widetilde{\nabla_y(w_i+y_i)}=\widetilde{\nabla_y w_i}(y)+e_i\chi_{Y^*},$$
and taking into account \cite[Corollary 2.10]{Cioranescu_Unfolding}, we have
\begin{eqnarray}\label{limit_gradient_Psi}
\widetilde{\nabla \Psi_{i\varepsilon}}\rightharpoonup
{1\over |Y|} \int_{Y^*}\left(e_i+\nabla_y w_i(y) \right)dy \quad \textrm{weakly in} \quad L^2(\Omega).
\end{eqnarray}
Due to that $w_i\in \mathbb{H}_{{\rm per}}/ \mathbb{R}$ (see \cite[Theorem 4.1]{Gahn}), let $\nabla_\Gamma \gamma_0(\Psi_{i\varepsilon})$ be the tangential gradient of $\gamma_0(\Psi_{i\varepsilon})$ on $\partial F_\varepsilon$ and  
we denote by $\mu_h^\varepsilon$ the above introduced linear form in the particular case in which $\displaystyle h\left({x\over \varepsilon}\right)=\nabla_\Gamma \gamma_0(\Psi_{i\varepsilon}(x))$.

From (\ref{definition_Psi}), we have
$$\nabla_\Gamma \gamma_0(\Psi_{i\varepsilon})=P_\Gamma e_i+\nabla_\Gamma w_i(y),$$
where $P_\Gamma e_i$ is defined on $\partial F$ and the tangential gradient of $w_i$ is given by
$$\nabla_\Gamma w_i:=P_\Gamma \nabla_y\tilde w_i=\nabla_y \tilde w_i-(\nabla_y \tilde w_i\cdot \nu)\nu \quad \text{on } \partial F,$$
where $\tilde w_i$ is an extension of $w_i$.

In this case, $\mu_{h}$ becomes 
$$\mu_h={1\over |Y|}\int_{\partial F}\left(P_\Gamma e_i+\nabla_\Gamma w_i(y)\right)d\sigma(y).$$
Then, using (\ref{convergence_mu}), we obtain
\begin{eqnarray}\label{limit_gradient_Beltrami_Psi}
\varepsilon \int_{\partial F_\varepsilon}\nabla_\Gamma \gamma_0(\Psi_{i\varepsilon}(x)) \varphi(x)d\sigma(x)=\langle \mu_{h}^\varepsilon,\varphi \rangle \to \langle \mu_{h},\varphi \rangle=\mu_h\int_{\Omega}\varphi(x)dx,\quad \forall \varphi\in W_0^{1,s}(\Omega).
\end{eqnarray}
On the other hand, it is not difficult to see that $ \Psi_{i\varepsilon}$ satisfies
\begin{equation}\label{system_Psi}
\left\{
\begin{array}{l}
\displaystyle -{\rm div}\left(\nabla \Psi_{i\varepsilon}\right)=0,   \text{\ in }\Omega_\varepsilon,\\[2ex]
\nabla \Psi_{i\varepsilon}\cdot \nu =\varepsilon\,\delta\,{\rm div}_{\Gamma}\left(\nabla_\Gamma \Psi_{i\varepsilon} \right),  \text{\ on }\partial F_\varepsilon.\end{array}
\right.
\end{equation}
Let $v\in \mathcal{D}(\Omega)$. Multiplying the first equation in (\ref{system_Psi}) by $vu_\varepsilon$, integrating by parts over $\Omega_\varepsilon$ and taking into account (\ref{test_function_boundary}), we get
\begin{eqnarray}\label{igualdad_problema_Psi}
&-&\varepsilon\,\delta\int_{\partial F_\varepsilon}\nabla_\Gamma \gamma_0(\Psi_{i\varepsilon})\cdot \nabla v\,\gamma_0(u_\varepsilon) d\sigma(x)-\varepsilon\,\delta\int_{\partial F_\varepsilon}\nabla _\Gamma\gamma_0(\Psi_{i\varepsilon})\cdot \nabla_\Gamma \gamma_0(u_\varepsilon) vd\sigma(x)\\ \nonumber
&=&\int_{\Omega_\varepsilon}\nabla \Psi_{i\varepsilon}\cdot \nabla v\,u_\varepsilon dx+\int_{\Omega_\varepsilon}\nabla \Psi_{i\varepsilon}\cdot \nabla u_\varepsilon vdx.
\end{eqnarray}
On the other hand, we multiply system (\ref{PDE}) by the test function $v \Psi_{i\varepsilon}$, integrating by parts over $\Omega_\varepsilon$ and taking into account (\ref{test_function_boundary}), we obtain
\begin{eqnarray*}
\dfrac{d}{dt}\left(\int_{\Omega}\chi_{\Omega_\varepsilon}\tilde u_\varepsilon v\tilde\Psi_{i\varepsilon}dx\right)+\varepsilon\,\dfrac{d}{dt}\left(\int_{\partial F_\varepsilon}
\gamma_{0}(u_\varepsilon)v\gamma_0(\Psi_{i\varepsilon})d\sigma(x)\right)+\int_{\Omega_\varepsilon} \nabla u_\varepsilon \cdot \nabla v\Psi_{i\varepsilon}dx+\int_{\Omega_\varepsilon} \nabla u_\varepsilon \cdot \nabla \Psi_{i\varepsilon}vdx
\\[2ex]
+\kappa \int_{\Omega}\chi_{\Omega_\varepsilon}\tilde u_\varepsilon v\tilde \Psi_{i\varepsilon} dx+\varepsilon\, \delta\int_{\partial F_\varepsilon} \nabla_\Gamma \gamma_0(u_\varepsilon)\cdot\nabla v\,\gamma_0(\Psi_{i\varepsilon})d\sigma(x)+\varepsilon\, \delta\int_{\partial F_\varepsilon} \nabla_\Gamma \gamma_0(u_\varepsilon)\cdot\nabla_\Gamma \gamma_0(\Psi_{i\varepsilon})vd\sigma(x)
\\[2ex]
 +\varepsilon\, \int_{\partial F_\varepsilon}
g(\gamma_{0}(u_\varepsilon))v\gamma_0(\Psi_{i\varepsilon})d\sigma(x)
 =0, 
 \end{eqnarray*}
 in $\mathcal{D}'(0,T)$. 
 
Using (\ref{igualdad_problema_Psi}), we have
 \begin{eqnarray*}
\dfrac{d}{dt}\left(\int_{\Omega}\chi_{\Omega_\varepsilon}\tilde u_\varepsilon v\tilde\Psi_{i\varepsilon}dx\right)+\varepsilon\,\dfrac{d}{dt}\left(\int_{\partial F_\varepsilon}
\gamma_{0}(u_\varepsilon)v\gamma_0(\Psi_{i\varepsilon})d\sigma(x)\right)+\int_{\Omega_\varepsilon} \nabla u_\varepsilon \cdot \nabla v\Psi_{i\varepsilon}dx-\int_{\Omega_\varepsilon}\nabla \Psi_{i\varepsilon}\cdot \nabla v\,u_\varepsilon dx
\\[2ex]
-\varepsilon\,\delta\int_{\partial F_\varepsilon}\nabla_\Gamma \gamma_0(\Psi_{i\varepsilon})\cdot \nabla v\,\gamma_0(u_\varepsilon) d\sigma(x)+\kappa \int_{\Omega}\chi_{\Omega_\varepsilon}\tilde u_\varepsilon v\tilde\Psi_{i\varepsilon} dx
+\varepsilon\, \delta\int_{\partial F_\varepsilon} \nabla_\Gamma \gamma_0(u_\varepsilon)\cdot\nabla v\,\gamma_0(\Psi_{i\varepsilon})d\sigma(x)
\\[2ex]
+\varepsilon\, \int_{\partial F_\varepsilon}
g(\gamma_{0}(u_\varepsilon))v\gamma_0(\Psi_{i\varepsilon})d\sigma(x)
 =0, 
 \end{eqnarray*}
 in $\mathcal{D}'(0,T)$.
 
We consider $\varphi\in C_c^1([0,T])$ such that $\varphi(T)=0$ and $\varphi(0)\ne 0$. Multiplying by $\varphi$ and integrating between $0$ and $T$, we have
 \begin{eqnarray}\label{New_formula}
 -\varphi(0)\left(\int_{\Omega}\chi_{\Omega_\varepsilon}\tilde u_\varepsilon(0)v\tilde\Psi_{i\varepsilon}dx\right)-\int_0^T\dfrac{d}{dt}\varphi(t)\left(\int_{\Omega}\chi_{\Omega_\varepsilon}\tilde u_\varepsilon(t)v\tilde\Psi_{i\varepsilon}dx\right)dt\\[2ex] \nonumber
-\varepsilon\varphi(0)\left(\int_{\partial F_\varepsilon}
\gamma_{0}(u_\varepsilon(0))v\gamma_0(\Psi_{i\varepsilon})d\sigma(x)\right)
-\varepsilon\int_0^T \dfrac{d}{dt}\varphi(t)\left(\int_{\partial F_\varepsilon}
\gamma_{0}(u_\varepsilon(t))v\gamma_0(\Psi_{i\varepsilon})d\sigma(x)\right)dt\\[2ex] \nonumber
+\int_0^T \varphi(t)\int_{\Omega} \tilde \xi_\varepsilon \cdot \nabla v\tilde \Psi_{i\varepsilon}dxdt-\int_0^T\varphi(t)\int_{\Omega}\widetilde{\nabla \Psi_{i\varepsilon}}\cdot \nabla v\,\tilde u_\varepsilon dxdt
\\[2ex]\nonumber
-\varepsilon\,\delta\int_0^T \varphi(t)\int_{\partial F_\varepsilon}\nabla_\Gamma \gamma_0(\Psi_{i\varepsilon})\cdot \nabla v\,\gamma_0(u_\varepsilon) d\sigma(x)dt+\kappa\int_0^T \varphi(t) \int_{\Omega}\chi_{\Omega_\varepsilon}\tilde u_\varepsilon v\tilde\Psi_{i\varepsilon} dxdt
\\[2ex] \nonumber
+\varepsilon\, \delta\int_0^T \varphi(t)\int_{\partial F_\varepsilon} \nabla_\Gamma \gamma_0(u_\varepsilon)\cdot\nabla v\,\gamma_0(\Psi_{i\varepsilon})d\sigma(x)dt
\\[2ex]
+\varepsilon\, \int_0^T \varphi(t)\int_{\partial F_\varepsilon}
g(\gamma_{0}(u_\varepsilon))v\gamma_0(\Psi_{i\varepsilon})d\sigma(x)dt
 =0.\nonumber
 \end{eqnarray}
Now, we have to pass to the limit, as $\varepsilon\to 0$. We will focus on the terms which involve the gradient and the tangential gradient. Taking into account (\ref{convergence_Psi_fuerte}), we reason as in \cite[Theorem 6.1]{Anguiano} for the others terms. 
 
 Firstly, using (\ref{convergence_gradiente2}), (\ref{convergence_Psi_fuerte}) and Lebesgue's Dominated Convergence Theorem, we have 
\begin{eqnarray*}\label{limite_1}
\int_0^T \varphi(t)\int_{\Omega} \tilde \xi_\varepsilon \cdot \nabla v\tilde \Psi_{i\varepsilon}dxdt\to \int_0^T \varphi(t)\int_{\Omega}\xi \cdot \nabla v\, x_idxdt,
\end{eqnarray*}
and by (\ref{converge_initial_data}), (\ref{limit_gradient_Psi}) and Lebesgue's Dominated Convergence Theorem, we obtain
\begin{eqnarray*}\label{limite_2}
\int_0^T\varphi(t)\int_{\Omega}\widetilde{\nabla \Psi_{i\varepsilon}}\cdot \nabla v\,\tilde u_\varepsilon dxdt\to {1\over |Y|}\int_0^T \varphi(t)\int_{\Omega}\left( \int_{Y^*}\left(e_i+\nabla_y w_i \right)dy\right)\cdot \nabla v\, udxdt.
\end{eqnarray*}
On the other hand, using (\ref{continuity1}) and (\ref{limit_gradient_Beltrami_Psi}), we can deduce
$$
\varepsilon\,\delta\int_{\partial F_\varepsilon}\nabla_\Gamma \gamma_0(\Psi_{i\varepsilon})\cdot \nabla v\,\gamma_0(u_\varepsilon) d\sigma(x)\to {\delta \over |Y|} \int_{\Omega}\left(\int_{\partial F}\left(P_\Gamma e_i+\nabla_\Gamma w_i \right)d\sigma(y)\right)\cdot\nabla v\, udx,
$$
which integrating in time and by Lebesgue's Dominated Convergence Theorem, we obtain
\begin{eqnarray*}\label{limite_3}
\varepsilon\,\delta\int_0^T \!\!\!\!\varphi(t)\int_{\partial F_\varepsilon}\nabla_\Gamma \gamma_0(\Psi_{i\varepsilon})\cdot \nabla v\,\gamma_0(u_\varepsilon) d\sigma(x)dt\to {\delta \over |Y|} \int_0^T \!\!\!\!\varphi(t)\int_{\Omega}\left(\int_{\partial F}\left(P_\Gamma e_i+\nabla_\Gamma w_i \right)d\sigma(y)\right)\cdot\nabla v\, udxdt.
\end{eqnarray*}
Similarly to the proof of (\ref{limit_Beltrami}) together with (\ref{convergence_Psi_fuerte}), we have
\begin{eqnarray*}\label{limite_4}
\varepsilon\, \delta\int_0^T \varphi(t)\int_{\partial F_\varepsilon} \nabla_\Gamma \gamma_0(u_\varepsilon)\cdot\nabla v\,\gamma_0(\Psi_{i\varepsilon})d\sigma(x)dt\to \delta {|\partial F|\over |Y|}\int_0^T \varphi(t)\int_{\Omega}\xi\cdot \nabla v\, x_i dxdt.
\end{eqnarray*}
Therefore, when we pass to the limit in (\ref{New_formula}), we obtain
\begin{eqnarray*}
 -\varphi(0)\left({|Y^*|\over |Y|}+{|\partial F|\over |Y|}\right)\left(\int_{\Omega} u(0)vx_i dx\right)-\left({|Y^*|\over |Y|}+{|\partial F|\over |Y|}\right)\int_0^T\dfrac{d}{dt}\varphi(t)\left(\int_{\Omega} u(t)v x_idx\right)dt\\[2ex] \nonumber
+\int_0^T \varphi(t)\int_{\Omega}\xi \cdot \nabla v\, x_idxdt-{1\over |Y|}\int_0^T \varphi(t)\int_{\Omega}\left( \int_{Y^*}\left(e_i+\nabla_y w_i \right)dy\right)\cdot \nabla v\, udxdt
\\[2ex]\nonumber
-{\delta \over |Y|} \int_0^T \!\!\!\!\varphi(t)\int_{\Omega}\left(\int_{\partial F}\left(P_\Gamma e_i+\nabla_\Gamma w_i \right)d\sigma(y)\right)\cdot\nabla v\, udxdt+\kappa{|Y^*|\over |Y|}\int_0^T \varphi(t) \int_{\Omega} u vx_idxdt
\\[2ex]
+\delta {|\partial F|\over |Y|}\int_0^T \varphi(t)\int_{\Omega}\xi\cdot \nabla v\, x_i dxdt
+{|\partial F|\over |Y|} \int_0^T \varphi(t)\int_{\Omega}
g(u(t))vx_idxdt
 =0.
 \end{eqnarray*}
Using Green's formula and equation (\ref{equation_xi}), we have
\begin{eqnarray*}
-\int_0^T \varphi(t)\int_{\Omega}\xi \cdot \nabla x_i\, vdxdt+{1\over |Y|}\int_0^T \varphi(t)\int_{\Omega}\left( \int_{Y^*}\left(e_i+\nabla_y w_i \right)dy\right)\cdot \nabla u\,vdxdt
\\[2ex]
+{\delta \over |Y|} \int_0^T \!\!\!\!\varphi(t)\int_{\Omega}\left(\int_{\partial F}\left(P_\Gamma e_i+\nabla_\Gamma w_i \right)d\sigma(y)\right)\cdot\nabla u\, vdxdt
-\delta {|\partial F|\over |Y|}\int_0^T \varphi(t)\int_{\Omega}\xi\cdot \nabla x_i\, vdxdt =0.
\end{eqnarray*}
The above equality holds true for any $v\in \mathcal{D}(\Omega)$ and $\varphi\in C_c^1([0,T])$. This implies that
 \begin{eqnarray*}
 -\left(1+\delta {|\partial F|\over |Y|}\right)\xi \cdot \nabla x_i +{1\over |Y|}\left( \int_{Y^*}\left(e_i+\nabla_y w_i \right)dy\right)\cdot \nabla u
+{\delta \over |Y|} \left(\int_{\partial F}\left(P_\Gamma e_i+\nabla_\Gamma w_i \right)d\sigma(y)\right)\cdot\nabla u =0,
\end{eqnarray*}
in $\Omega \times (0,T).$ We conclude that
 \begin{eqnarray}\label{identificacion_xi}
\left(1+\delta {|\partial F|\over |Y|}\right) {\rm div}\xi= {\rm div}\left(Q \nabla u\right),
\end{eqnarray}
where $Q=((q_{ij}))$, $1\leq i,j\leq N$, is given by 
\begin{eqnarray*}\label{A_limite}
q_{ij}={1\over |Y|}\left( \int_{Y^*}\left(e_{i}+\nabla_y w_i \right)\cdot e_j \,dy
+\delta  \int_{\partial F}\left(P_\Gamma e_{i}+\nabla_\Gamma w_i \right)\cdot P_\Gamma e_j \,d\sigma(y)\right).
\end{eqnarray*}
Observe that if we multiply system (\ref{system_eta}) by the test function $w_j$, integrating by parts over $Y^*$, we obtain
\begin{eqnarray*}
\int_{Y^*}(e_i+\nabla_y w_i)\cdot \nabla_y w_j dy+\delta \int_{\partial F}(P_\Gamma e_i+\nabla_\Gamma w_i)\cdot \nabla_\Gamma w_j d\sigma(y)=0,
\end{eqnarray*}
then we conclude that $q_{ij}$ is given by (\ref{matrix}).

{\bf Step 3.} Finally, thanks to (\ref{equation_xi}) and (\ref{identificacion_xi}), we observe that $u$ satisfies the first equation in (\ref{limit_problem}). A weak solution of (\ref{limit_problem}) is any function $u$, satisfying
\begin{equation*}
u\in {C}([0,T];L^{2}\left(  \Omega\right)  ),\quad \text{for all }T>0,
\end{equation*}
\begin{equation*}
u\in L^{2}(0,T;H_{0}^{1}\left(  \Omega\right)  )\cap L^{q}(0
,T;L^{q}\left(  \Omega\right)  ),\quad \text{for all }T>0,
\end{equation*}
\begin{equation*}
\displaystyle \left({|Y^*|\over |Y|}+{|\partial F| \over |Y|} \right) \dfrac{d}{dt}(u(t),v)+(Q\nabla u(t),\nabla v)+{|Y^*|\over |Y|}\kappa(u(t),v)+{|\partial F| \over |Y|}(g(u(t)),v)=0,\quad \text{in  }\mathcal{D}'(0,T),
\end{equation*}
for all $v\in H_0^1(\Omega)\cap L^q(\Omega)$, and
\begin{equation*}
 u(0)=u_0.
\end{equation*}
Due to that the homogenized matrix $Q$ is positive-definite (see \cite[Theorem 4.1]{Gahn}), applying a slight modification of \cite[Chapter 2,Theorem 1.4]{Lions}, we obtain that the problem (\ref{limit_problem}) has a unique solution, and therefore Theorem \ref{Main} is proved.

\begin{remark}
It is worth remaking that if we consider a nonlinear term $f(u_\varepsilon)$ in the first equation in (\ref{PDE}) which satisfies the same assumptions as $g$, we obtain Theorem \ref{Main} with an additional term $\displaystyle{|Y^*|\over |Y|}f(u)$ in the first equation in (\ref{limit_problem}).
\end{remark}


\begin{thebibliography}{10}

\bibitem{Amar_Gianni} M. Amar, R. Gianni, Laplace-Beltrami operator for the heat conduction in polymer coating of electronic devices, Discrete and Continuous Dynamical Systems Series B, 23, No. 4 (2018) 1739-1756.

\bibitem{Amar_Gianni2} M. Amar, R. Gianni, Error estimate for a homogenization problem involving the Laplace-Beltrami operator, Mathematics and Mechanics of Complex Systems, 6, No. 1 (2018) 41-59.

\bibitem{Anguiano} M. Anguiano, Existence, uniqueness and homogenization of nonlinear parabolic problems with dynamical boundary conditions in perforated media, Mediterr. J. Math. (2020) 17:18.

\bibitem{Cioranescu_Unfolding} D. Cioranescu, A. Damlamian, G. Griso, The periodic Unfolding Method in Homogenization, SIAM Journal on Mathematical Analysis, 40, No. 4 (2008) 1585-1620.

\bibitem{Ciora2} D. Cioranescu, P. Donato, Homog\'en\'eisation du probl\`eme de Neumann non homog\`ene dans des ouverts perfores, Asymptotic Analysis, 1, (1988) 115-138.
%
\bibitem{CioraDonato} D. Cioranescu, P. Donato, An Introduction to Homogenization, Oxford Lectures Series in Mathematics and its Applications, 17, New York, 1999.

\bibitem{Cioranescu} D. Cioranescu, J. Saint Jean Paulin, Homogenization in open sets with holes, J. Math. Anal. Appl., 71, (1979) 590-607.

\bibitem{Conca} C. Conca, J.I. D\'iaz, A. Li\~n\'an, C. Timofte, Homogenization in chemical reactive flows, Electronic Journal of Differential Equations. 40 (2004), 1-22.

\bibitem{Gahn} M. Gahn, Multi-scale modeling of processes in porous media-coupling reaction-diffusion processes in the solid and the fluid phase and on the separating interfaces, Discrete and Continuous Dynamical Systems Series B, 24, No. 12 (2019) 6511-6531.

\bibitem{Gal} C.G. Gal, The role of surface diffusion in dynamic boundary conditions: where do we stand?, Milan J. Math. 83, (2015) 237-278.

\bibitem{Gal2} C.G. Gal, J. Shomberg, Coleman-Gurtin type equations with dynamic boundary conditions, Phys. D. 292, (2015) 29-45.


\bibitem{Goldstein} G.R. Goldstein, Derivation and physical interpretation of general boundary conditions, Adv. Differential Equations, 11, No. 4 (2006) 457-480.

\bibitem{Graf_Peter} I. Graf, M.A. Peter, Diffusion on surfaces and the boundary periodic unfolding operator with an application to carcinogenesis in human cells, SIAM J. Math. Analysis, 46, No. 4 (2014) 3025-3049.

\bibitem {Lions}J.L. Lions, Quelques M\'{e}thodes de R\'{e}solution des
Probl\`{e}mes aux Limites Non lin\`{e}aires, Dunod, 1969.

\bibitem{Lions_Magenes} J.L. Lions, E. Magenes, Non-Homogeneous Boundary Value Problems and Applications, Vol. 1, Springer-Verlag Berlin Heidelberg, 1972.

\bibitem{Robinson} J.C. Robinson, Infinite-dimensional dynamical systems. Cambridge University Press, 2001.

\bibitem{Soko_Zole} J. Sokolowski, J.-P. Zolesio, Introduction to Shape Optimization, Springer-Verlag, 1992.

\bibitem{Strichartz} R.S. Strichartz, Analysis of the Laplacian on the complete Riemannian manifold, Journal of Functional Analysis, 52, (1983) 48-79.

\bibitem{Tartar} L. Tartar, Problèmes d'homog\'en\'eisation dans les \'equations aux d\'eriv\'ees partielles, Cours Peccot Collège de France, 1977.

\bibitem{Vanni} M. Vanninathan, Homogenization of eigenvalues problems in perforated domains, Proc. Indian Acad. of Sciences, 90, (1981) 239-271.

\end{thebibliography}
\end{document}